\documentclass[11pt]{article}
\usepackage{amsthm,amsmath,amssymb}

{\catcode `\@=11 \global\let\AddToReset=\@addtoreset}
\AddToReset{equation}{section}

\usepackage{graphicx}
\usepackage{epsfig}
\usepackage{color}

\usepackage[latin1]{inputenc}

\newtheorem{theorem}{Theorem}
\newtheorem{proposition}{Proposition}
\newtheorem{@definition}{Definition}
\newenvironment{definition}{\begin{@definition}\rm}{\end{@definition}}
\newtheorem{@remark}{Remark}
\newenvironment{remark}{\begin{@remark}\rm}{\end{@remark}}

\def\R{\mathbb{R}}
\def\Z{\mathbb{Z}}
\def\i{\mathrm i}
\def\d{\mathrm d}

\def\I{\mathbb{I}}

\def\e{\mathrm e}
\def\E{\mathrm E}

\def\1{{\bf 1}}
\def\0{{\bf 0}}

\def\vep{\varepsilon}

\newcommand{\nn}{\nonumber}
\newcommand{\noi}{\noindent}

\newcommand{\mbt}{\boldsymbol{t}}

\newcommand{\mbs}{\boldsymbol{s}}

\newcommand{\mbu}{\boldsymbol{u}}

\newcommand{\mbx}{\boldsymbol{x}}

\newcommand{\mby}{\boldsymbol{y}}

\def\limd{\renewcommand{\arraystretch}{0.5}
\begin{array}[t]{c}
\stackrel{\rm d}{\longrightarrow} \\
\end{array}\renewcommand{\arraystretch}{1}}

\def\limfdd{\renewcommand{\arraystretch}{0.5}
\begin{array}[t]{c}
\stackrel{\rm fdd}{\longrightarrow} \\
\end{array}\renewcommand{\arraystretch}{1}}

\def\eqfdd{\renewcommand{\arraystretch}{0.5}
\begin{array}[t]{c}
\stackrel{\rm fdd}{=} \\
\end{array}\renewcommand{\arraystretch}{1}}

\def\neqfdd{\renewcommand{\arraystretch}{0.5}
\begin{array}[t]{c}
\stackrel{\rm fdd}{\neq} \\
\end{array}\renewcommand{\arraystretch}{1}}

\begin{document}


\title{Scaling transition for singular linear random fields on $\Z^2$: \\
spectral approach }
\author{Donatas Surgailis}
\date{\today
\\  \small
\vskip.2cm
Vilnius University, Faculty of Mathematics and Informatics,
Naugarduko 24, 03225  Vilnius, Lithuania \\
}
\maketitle

\begin{abstract}

\medskip

We study partial sums limits of linear  random fields $X$ on $\Z^2 $ with spectral density $f(\mbx)
$ tending to $\infty,\, 0$ or to both (along different subsequences) as $\mbx \to (0,0)$.
The above behaviors are termed
(spectrum) long-range dependence, negative dependence, and long-range negative dependence, respectively, and
assume an  anisotropic power-law form of  $f(\mbx)$   near the origin. The partial sums are taken  over rectangles whose
sides increase as $\lambda$ and  $\lambda^\gamma $, for any fixed $\gamma >0$. We prove that for above $X$  the partial sums or  scaling  limits exist for any $\gamma>0$ and exhibit a scaling transition  
at some  $\gamma = \gamma_0>0$; moreover,  the `unbalanced' scaling limits ($\gamma\ne\gamma_0$) are  Fractional Brownian  Sheet
with Hurst parameters taking values from $[0,1]$. The paper extends \cite{ps2015, pils2017, sur2020}  
to the above spectrum dependence conditions and/or more general values of Hurst parameters.

\medskip

{\small

\noi {\it Keywords:} Linear random field, Gaussian random field, spectral density,  long-range dependence, negative dependence, long-range  negative dependence,
hyperbolic dependence,
self-similarity,  anisotropic scaling, scaling transition, fractional Brownian
sheet,

}

%

\end{abstract}

\section{Introduction}

A stationary random field (RF) $ X = \{ X(\mbt); \mbt \in \Z^d \} $ with finite second moment  and covariance $r(\mbt)
= {\rm Cov}(X(\0), X(\mbt)) $ is said
(i) long-range dependent (LRD)  if $\sum_{\mbt\in \Z^d}  |r(\mbt)|  = \infty$; \ (ii) short-range dependent (SRD)
if $\sum_{\mbt\in \Z^d}  |r(\mbt)|  < \infty, \
$  $ \sum_{\mbt\in \Z^d}  r(\mbt) \ne 0$, and  (iii) negatively dependent (ND)
if  $\sum_{\mbt\in \Z^d}  |r(\mbt)| < \infty, \ \sum_{\mbt\in \Z^d}  r(\mbt)  = 0 $
(for RF indexed by continuous argument, the above concepts are analogously defined with sums replaced by integrals over $\R^d $).
The above classification, albeit not completely satisfactory, is important in limit theorems, see \cite{gir2012, lah2016, sur2020}. Related
but not equivalent definitions of LRD, SRD and ND are given in terms of spectrum or moving-average  coefficients of RF. In the sequel we refer
to the above covariance or moving-average characterizations as {\it spatial domain} properties, to be distinguished from {\it frequency  domain (spectrum) } characterizations used in the literature. A very general form of such frequency domain
concepts assumes  the existence of spectral density $f(\mbu), \mbu \in \Pi^d
:= [-\pi,  \pi]^d $ which is bounded outside of the origin and
such that  (i') $\lim_{\mbu \to \0} f(\mbu) = \infty $ (spectrum LRD),
(ii') 
$\lim_{\mbu \to \0} f(\mbu) >0  $  (spectrum SRD), and
(iii') 
$\lim_{\mbu \to \0} f(\mbu) = 0  $ (spectrum ND).  Moreover, in cases (i') and (iii') it is often additionally
assumed that spectral density varies regularly as $\mbu \to 0$ \cite{gir2012, pipi2017, samo2016}. It is clear that in dimensions $d \ge 2 $
a function can increase regularly to $\infty $ along one direction and vanish along another one. Such a behavior is not
reflected in (i')-(iii'),  calling for a new concept
(iv')  $0 = \liminf_{\mbu \to \0} f(\mbu) < \limsup_{\mbu \to \0} f(\mbu)  =  \infty $  which we term
{\it spectrum long-range negative dependence }
(spectrum LRND).
Properties (i'), (iii') and (iv') can be jointly dubbed as    {\it singular } behaviors, the singularity limited to the origin and including both infinite and/or zero limits.

The present paper studies anisotropic partial sums limits 
for a class of stationary linear RFs  $X$ on $\Z^2 $ with spectrum
LRD/ND/LRND properties. The partial  sums
\begin{equation}\label{Sgamma}
S_{\lambda,\gamma}(\mbx)  :=  \
\sum_{\mbt \in K_{[\lambda x_1,\lambda^{\gamma} x_2]}}
X(\mbt), \qquad \mbx =  (x_1, x_2) \in \R^2_+
\end{equation}
are taken over rectangles $K_{[\lambda x_1, \lambda^{\gamma}x_2]} := \{ \mbt  \in \Z^2: 1\le t_1 \le \lambda x_1, 1\le t_2 \le \lambda^{\gamma}x_2 \}
\subset \Z^2 $
whose sides increase as $O(\lambda) $ and $O(\lambda^\gamma)$ as
$\lambda \to \infty$, for a fixed
$\gamma >0$.  Our main object
are anisotropic scaling limits
\begin{equation} \label{VgammaD}
d^{-1}_{\lambda,\gamma} S_{\lambda,\gamma}(\mbx)
 \ \limfdd \ V_\gamma (\mbx), \quad \mbx \in \R^2_+, \quad \lambda \to \infty,
\end{equation}
where  $d_{\lambda,\gamma} \to \infty $ is a normalization. It was proved in  \cite{ps2015, ps2016, pils2016, pils2017, bier2017, pils2020, sur2020} and other works that for a large class of RFs
nontrivial scaling limits  $V_\gamma (\mbx)$ in \eqref{VgammaD} exist for any $\gamma >0$. Moreover,
the limit family $\{V_\gamma; \gamma >0\}$ may exhibit a trichotomy called {\it scaling transition} at some (nonrandom) point $\gamma = \gamma_0 >0$ meaning
that $V_\gamma = V_\pm \ (\gamma \lessgtr \gamma_0)$  do not depend on $\gamma $ and are different; see Definition \ref{deftrans}. 
For  linear RFs,  these works can be summarized as saying that a large class of
 RFs on $\Z^2 $ parade two types of the limit behavior in  \eqref{VgammaD}: either (I) $V_\gamma $ is a standard Brownian Sheet $B_{1/2,1/2} $
for any $\gamma >0$, or  (II)  scaling transition at some $\gamma_0  >0$ exists such that all {\it unbalanced} limits
$V_\gamma, \gamma \ne \gamma_0 $ are Fractional Brownian  Sheet (FBS) $B_{H_1,H_2}  $ with Hurst parameters
$(H_1,H_2) = (H^+_1,H^+_2)\, (\gamma > \gamma_0),\  = (H^-_1,H^-_2) \, (\gamma < \gamma_0)$,
$ (H^+_1,H^+_2) \neq  (H^-_1,H^-_2) $. Behavior (I) is rather standard under  SRD conditions, see \cite{dam2017, dam2021, ps2016} and the references 
therein. Behavior (II) is more interesting
and was established under LRD and ND conditions in \cite{ps2015, sur2020}.  However, these results exclude LRND singularity and 
some particular (boundary)  values of parameters.

\begin{figure}[ht!]
\begin{center}
\begin{tabular}{cc}
\includegraphics[width=0.40\textwidth]{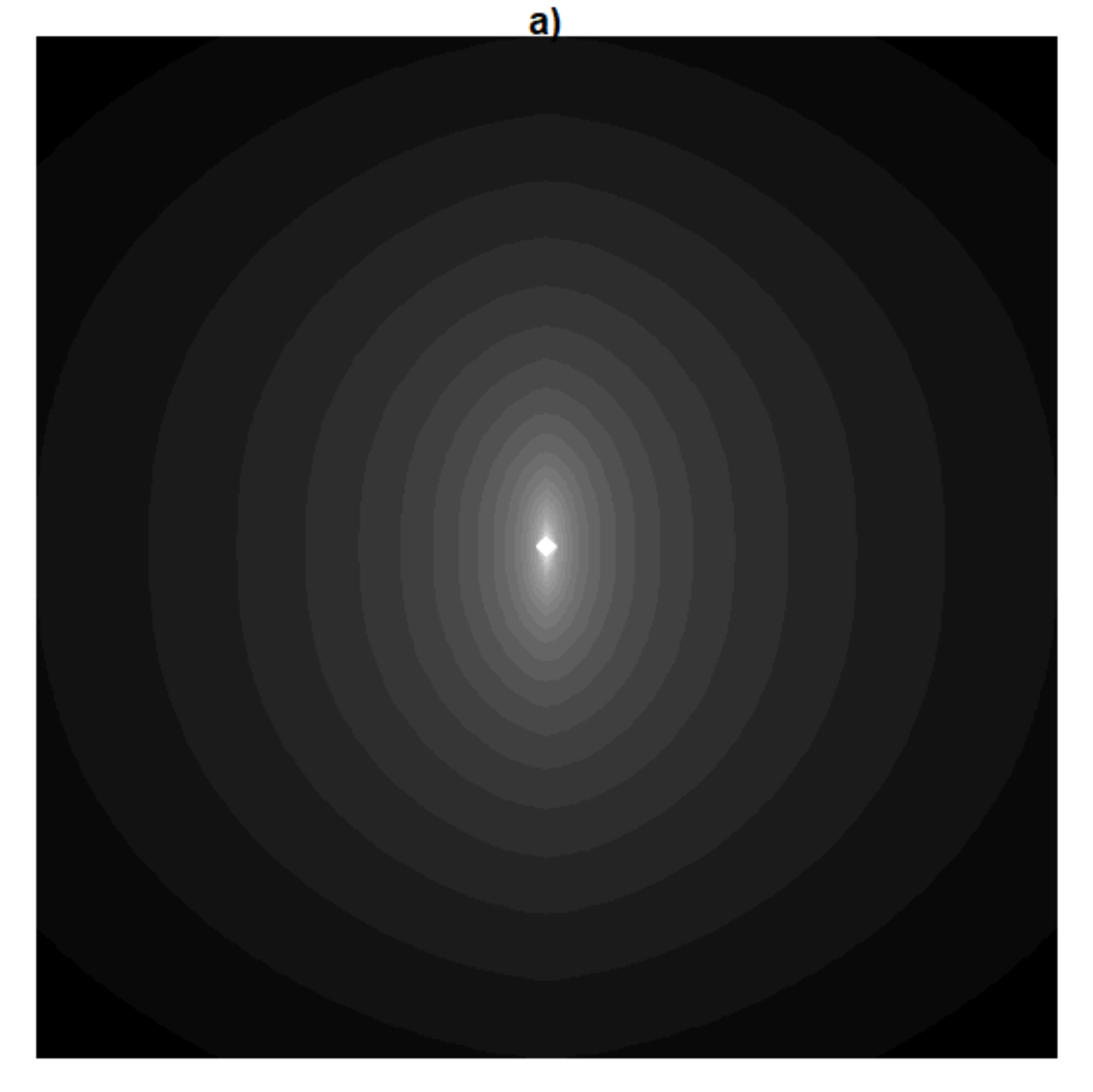}
&\hskip1cm \includegraphics[width=0.40\textwidth]{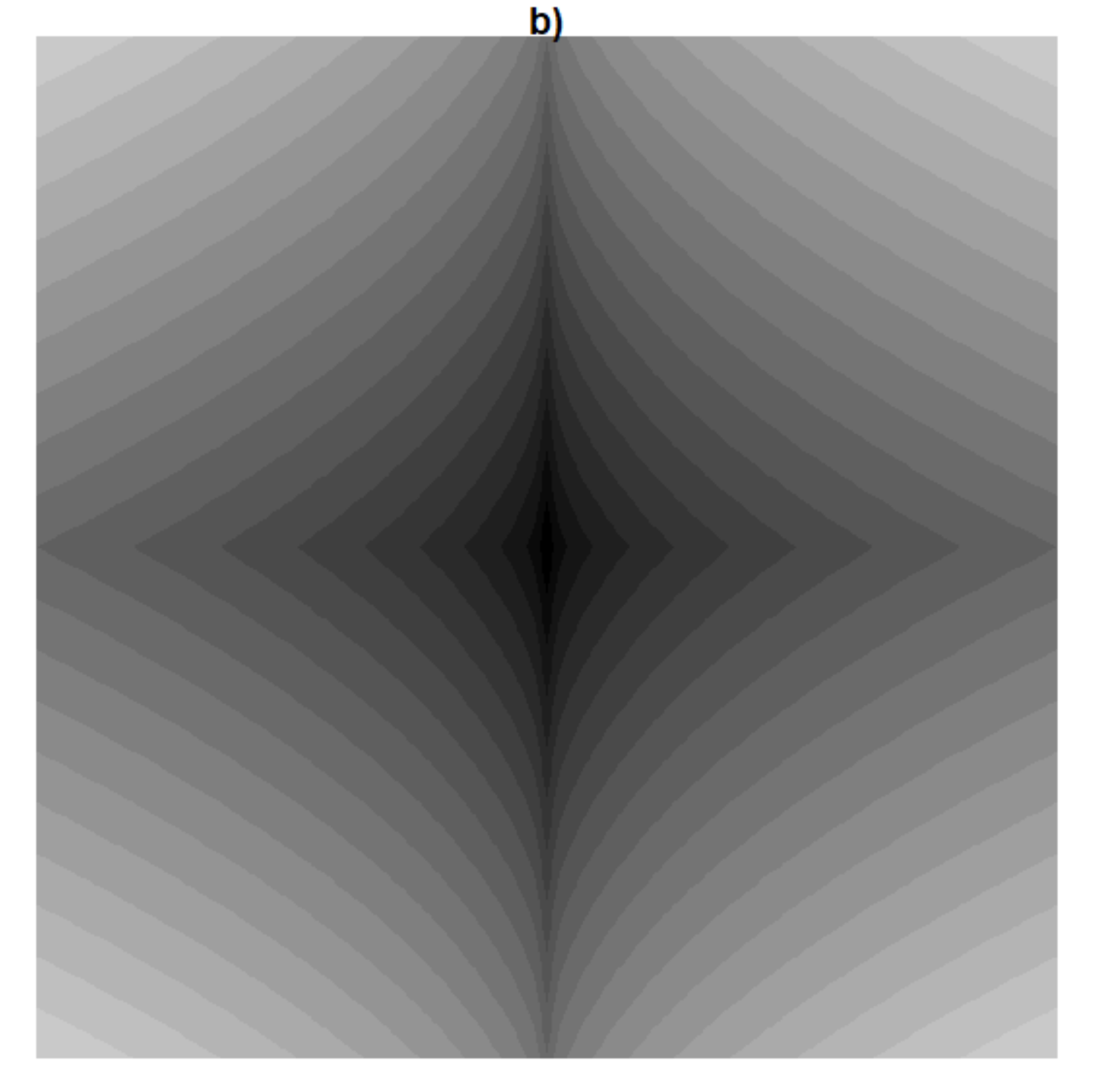}\\
 \includegraphics[width=0.40\textwidth]{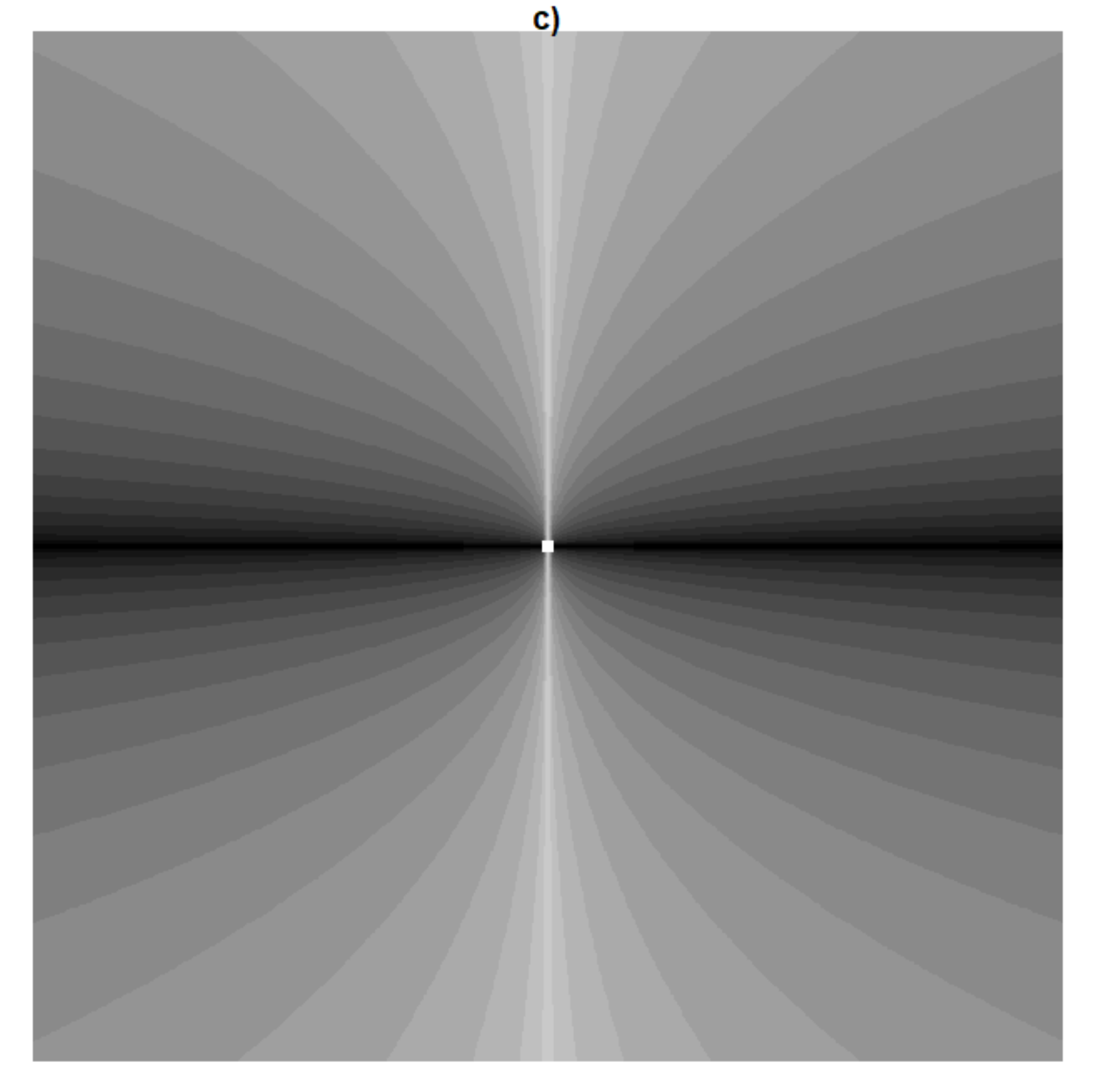}
&\hskip1cm \includegraphics[width=0.40\textwidth]{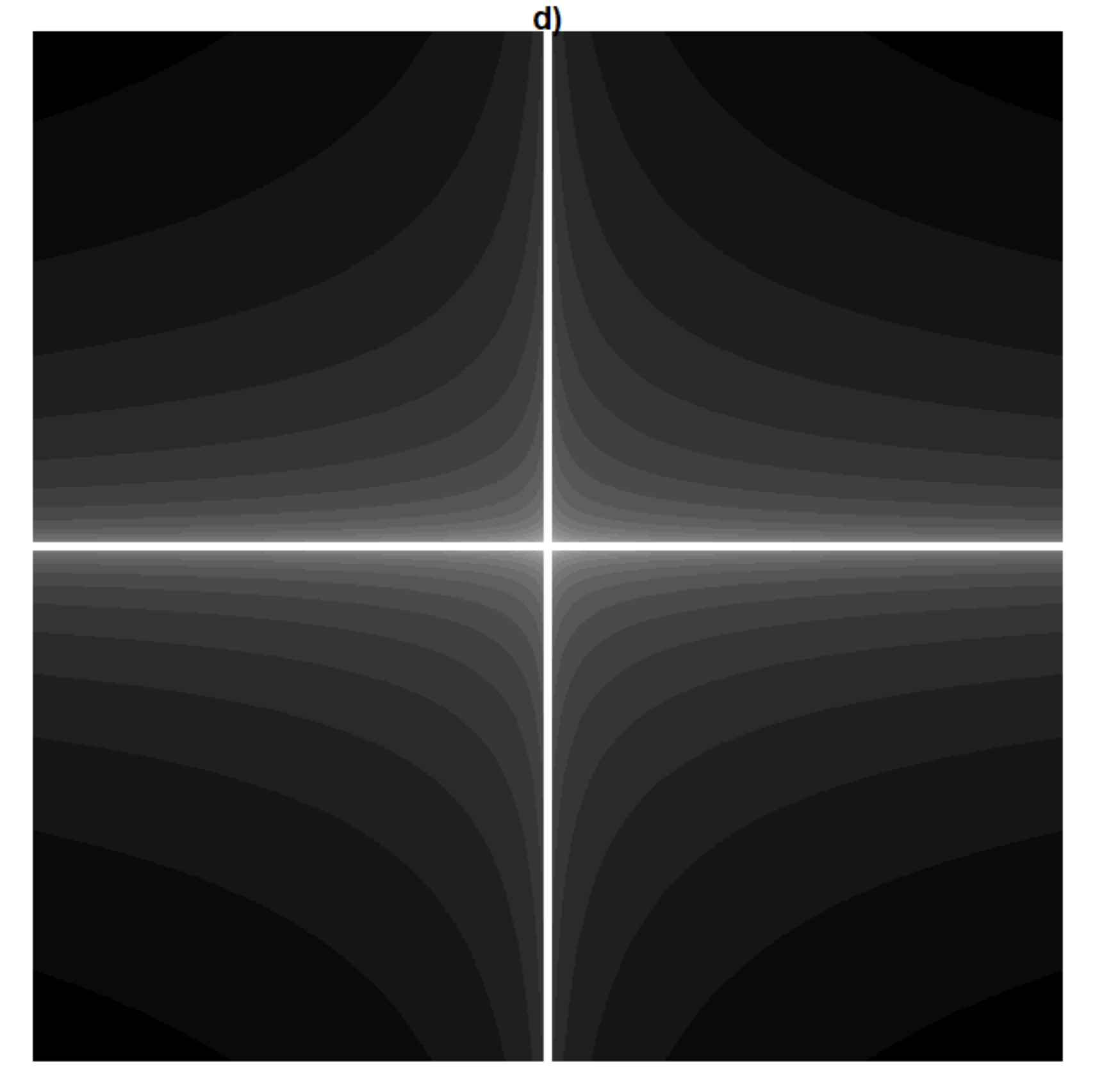}
\\
\end{tabular}
\end{center}

\noi {\small Figure 1. Level graphs of spectral density under LRD, ND, LRND and  `hyperbolic dependence' conditions.
a) $f(\mbx) = \rho^{-1}_2(\mbx)$, $\upsilon_1 = 0.8, \upsilon_2 = 1.2$; \
b) $f(\mbx) =  \rho(\mbx), \upsilon_1 = 0.5, \upsilon_2 = 1$;  \ c) $f(\mbx) =  \rho^{-1}(\mbx) |x_2/\rho(\mbx)^{1/\upsilon_2}|^{0.5},
\upsilon_1 = 0.5, \upsilon_2 = 1$; \ d) $f(\mbx) = |x_1|^{-0.2} |x_2|^{-.5}$
}

\end{figure}

The present paper extends \cite{ps2015, pils2017, sur2020} and some other work on type (II) behavior
for linear RFs $X$
under spectrum LRD,  ND, and LRND conditions. It appears that the spectral approach 
 used in this work 
is more efficient and simple compared to `time-domain' approach
in \cite{pils2017, sur2020} and other papers; moreover, it can be used to obtain the covariance structure of the scaling limits
in a non-Gaussian context as well. 
It is assumed that $X$ 
have a moving-average (MA) representation
\begin{equation}\label{Xlin}
X(\mbt) = \sum_{\mbs \in \Z^2 } a(\mbt- \mbs) \vep(\mbs), \qquad \mbt \in \Z^2,
\end{equation}
with standardized i.i.d. innovations $\{ \vep(\mbs); \mbs \in \Z^2\}, \E \vep(\mbs) = 0, \E \vep(\mbs)^2 =1 $ and
deterministic coefficients $a(\mbt), \mbt \in \Z^2. $ The latter coefficients satisfy the necessary condition $\sum_{\mbt \in \Z^2} a(\mbt)^2 < \infty $
for the convergence of \eqref{Xlin} but
all other conditions for various behaviors in \eqref{VgammaD} in this paper refer to the spectral density of $X$ written as the squared Fourier transform
$f(\mbx) = (2\pi)^2 |\widehat a(\mbx)|^2, $ 
$\widehat a(\mbx) := (2\pi)^{-2} \sum_{\mbt \in \Z^2} \e^{\i \mbt \cdot \mbx } a(\mbt), \ \mbx \in \Pi^2 $ of MA coeffiecients. 
For identification of the covariance structure of the limit Gaussian RF in  \eqref{VgammaD}, 
the basic assumption 
is that $f(\mbx) $ behaves like a power function
\begin{eqnarray}\label{rhosp}
\rho(\mbx) := |x_1|^{\upsilon_1} + |x_2|^{\upsilon_2}, \qquad \mbx = (x_1,x_2) \in \R^2
\end{eqnarray}
at the origin, where $\upsilon_i >0, i=1,2 $ are given parameters, or,
more precisely, 
\begin{equation}\label{fL}
f(\mbx) \sim  L(\mbx)/\rho(\mbx)  \quad \text{(`spectrum LRD')}  \quad \text{and}
\quad f(\mbx) \sim  L(\mbx) \rho(\mbx)  \quad   \text{(`spectrum ND')},
\end{equation}
where $L(\mbx) $ is an `angular function' satisfying some boundedness conditions. 
Clearly,
the first relation in \eqref{fL} implies spectrum LRD provided $L(\mbx)$ is bounded from below, while the second
relation  relation in \eqref{fL} implies spectrum ND provided $ L(\mbx)$ is bounded from above.
The above choice of 
$\rho(\mbx) $  as the `radial function' is rather flexible: we can use other forms of \eqref{rhosp}, e.g. with
\begin{eqnarray}\label{rhosp2}
\rho_p(\mbx) := (|x_1|^{p \upsilon_1} + |x_2|^{p \upsilon_2})^{1/p}, \quad p >0, \quad
\mbx = (x_1,x_2) \in \R^2
\end{eqnarray}
instead of \eqref{rhosp} as well, since the change of $\rho (\mbx)$ to  $\rho_p (\mbx)$  amounts
to a  change 
of the angular function which does not affect the limit results. 
Spectrum
LRND may arise if these boundedness conditions on $L(\mbx)$ are violated, e.g. when  $L(\mbx) $ decreases as $|x_i|^{\mu_i} $ with exponent
$\mu_i  >0$ when $x_i  \to 0$  in the first relation in \eqref{fL}.
We show that in this case the limits
in \eqref{VgammaD} may exist, together with a scaling transition and the limit Gaussian RF depends on parameters
$\upsilon_i, \mu_i, i=1,2$. We also consider the case when the spectral density asymptotically
factorizes into  a product of power functions $|x_i|^{-\upsilon_i}, i=1,2 $. 
This form of singularity
(studied in \cite{ps2015, dam2017} and termed  {\it hyperbolic dependence} in the present  paper)
allows for spectrum  LRD,  ND, or LRND depending on the sign of $\upsilon_i$,
but leads to very different limit results (the absence of scaling transition).   
Fig. 1 illustrates different types of singular spectral densities studied in this paper.
We also expect that our results can be extended to RFs on
$d$-dimensional lattice, $d \ge 3 $ arbitrary although the description of the limit RFs is more complicated, see
\cite{bier2017, sur2019}.

The rest of the paper is organized as follows. Sec.~2 provides rigorous assumptions on $f(\mbx) $
and some preliminary facts used in the subsequent text.  Sections 3-5 contain the main Theorems
\ref{thmGaussLRD},  \ref{thmGaussND}, \ref{thmGaussLRND} and \ref{farimaLM}, under respective LRD, ND, LRND and hyperbolic dependence premises.

\smallskip

{\it Notation.} In this paper, $\limfdd, \eqfdd, $ and $\neqfdd$,
denote respectively the weak convergence,  equality, and inequality, of finite dimensional distributions. $C$ stands for a generic positive constant which may assume different values at various locations and whose precise value has no importance.
$\R_+ := (0,\infty), \, \R^2_0 := \R^2 \setminus \{(0,0)\},  \Pi := [-\pi,\pi], \,  |\mbx| := |x_1| + |x_2|,
\mbx \cdot \mby := x_1 y_1 + x_2 y_2, \mbx = (x_1,x_2), \mby = (y_1,y_2), 
\1 := (1,1), \0 := (0,0)$.
$\I(A)$ denotes the indicator  function of a set $A$.

\section{Preliminaries and assumptions}

\begin{definition} \cite{ps2015, ps2016, pils2017} \label{deftrans} {\rm  We say that a stationary RF $ X = \{ X(\mbt); \mbt \in \Z^2 \} $ {\it exhibits  scaling transition} if non-trivial limits in \eqref{VgammaD} exist for any $\gamma >0$  and
there exists a
$\gamma_0 >0$  such that 
\begin{equation}\label{Vtrans}
V_\gamma \eqfdd V_+ \ \ (\gamma >\gamma_0), \qquad
V_\gamma \eqfdd V_- \ \ (0<\gamma < \gamma_0), \qquad V_+ \neqfdd a V_- \ \ (\forall \, a >0).
\end{equation}
In such a case, RFs $V_\pm $ will be called the {\it unbalanced} limits, and RF $V_0$ the {\it well-balanced} limit of $X$. }
\end{definition}

A weaker version of Definition \ref{deftrans} in \cite{dam2021}
requires that $V_\gamma $ is independent of $\gamma $ for $\gamma > \gamma_0$ and $\gamma < \gamma_0$ in a small neighborhood of $\gamma_0$.
There exist RFs which exhibit scaling transition in the above weaker sense at several (possibly, infinite) number of points $\gamma >0$ \cite{pils2016, dam2021}. However, for RFs discussed in this paper Definition  \ref{deftrans}  suffices.

Recall that the covariance function of stationary zero-mean  RF $X = \{ X(\mbt); \mbt \in \Z^2 \} $
is the Fourier transform of its spectral
density, viz., $r (\mbt)  := \E X(\0) X(\mbt) = \int_{\Pi^2} \e^{ \i \mbt \cdot \mbu} f(\mbu)  \d \mbu, \mbt \in \Z^2$. Whence, the covariance
function of the normalized partial sums RF in \eqref{Sgamma} writes as
\begin{eqnarray} \label{Rcov}
R_{\lambda, \gamma}(\mbx, \mby)
&:=&d^{-2}_{\lambda,\gamma}
\E S_{\lambda,\gamma}(\mbx) S_{\lambda, \gamma}(\mby)
\ =  \ d^{-2}_{\lambda,\gamma}
\int_{\Pi^2} \prod_{i=1}^2 D_{[\lambda^{\gamma_i} x_i]}(u_i) \overline{D_{[\lambda^{\gamma_i} y_i]}(u_i)} 
f(\mbu) \d \mbu, \quad \mbx, \mby \in \R^2_+,
\end{eqnarray}
where  $(\gamma_1, \gamma_2)  := (1,\gamma) $ and $D_n(u)
:= \sum_{t=1}^n  \e^{\i tu}  = (\e^{\i u}- \e^{\i (n+1)u})/(1- \e^{\i u}),  u \in \Pi$ (the real-valued function
$D_n(u)\e^{-\i (n+1) u/2} $
is the Dirichlet kernel).

\begin{proposition} Let $X$ be a linear RF in \eqref{Xlin} with  spectral  density $f$ and $S_{\lambda,\gamma}$ be the partial sum RF in
\eqref{Sgamma}, $\gamma >0$. Assume that 
\begin{equation} \label{Rlim}
R_{\lambda, \gamma}(\mbx, \mby)   \ \to \  R_\gamma (\mbx, \mby), 
\quad \lambda \to \infty \quad
(\forall \, \mbx, \mby \in \R^2_+)
\end{equation}
and
\begin{equation}\label{LindX}
|g_\lambda|_\infty :=
\sup_{\mbu \in \Z^2} \Big|\sum_{\mbt \in K_{[\lambda x_1,\lambda^{\gamma} x_2]}} a(\mbt - \mbu)\Big| = o(d_{\lambda,\gamma}), \qquad \lambda \to \infty.
\end{equation}
Then \eqref{VgammaD} holds, where $V_\gamma $ is a Gaussian RF on  $\R^2_+$ with zero mean $\E V_\gamma(\mbx) = 0$  and
covariance $\E V_\gamma(\mbx) V_\gamma (\mby) = R_\gamma (\mbx, \mby), \mbx, \mby \in \R^2_+$.

\end{proposition}

\noi {\bf Proof.} The finite-dimensional convergence in \eqref{VgammaD} is equivalent to  one-dimensional convergence
${\cal S}_{\lambda, \gamma} :=
d^{-1}_{\lambda,\gamma} \sum_{j=1}^m \theta_j $  $ S_{\lambda,\gamma}(\mbx_j) \limd \sum_{j=1}^m \theta_j V_{\gamma}(\mbx_j) =: {\cal V}_\gamma
$ for any $\theta_j \in  \R,  \mbx_j \in  \R^2_+,  j=1, \cdots, m, m \ge 1 $. From \eqref{Rlim} we have
$\E  ({\cal S}_{\lambda, \gamma})^2 \to \sum_{j,j'=1}^m \theta_j \theta_{j'} R_\gamma (\mbx_j, \mbx_{j'}) =:  {\cal R}_\gamma \ge 0$.
Let ${\cal R}_\gamma > 0$ then  ${\cal S}_{\lambda, \gamma} = d^{-1}_{\lambda,\gamma} \sum_{\mbu \in \Z^2} g_{\lambda} (\mbu) \vep(\mbu)$ is a normalized  weighted linear form in i.i.d. r.v.s with weights $g_\lambda (\mbu) := \sum_{j=1}^m \theta_j
\sum_{\mbt \in K_{[\lambda x_{1j},\lambda^{\gamma} x_{2j}]}} a(\mbt - \mbu), \mbx_j = (x_{1j}, x_{2j}) $. Relation
\eqref{LindX} implies the Lindeberg condition, viz.,  for any $\tau >0$
\begin{equation}\label{LindX1}
\sum_{\mbu \in \Z^2} g^2_{\lambda}(\mbt) \E  [\vep(\mbu)^2 \I( |g_{\lambda}(\mbu) \vep(\mbu)| >
\tau d_{\lambda,\gamma})] = o(d^2_{\lambda,\gamma}), \qquad \lambda \to \infty
\end{equation}
since $\E \vep(\mbu)^2 \I( |g_{\lambda} (\mbu) \vep(\mbu)| >
\tau d_{\lambda,\gamma}) \le
\E \vep(\mbu)^2 \I\big( |\vep(\mbu)| > \tau d_{\lambda,\gamma}/|g_\lambda|_\infty
\big) \to 0 $ 
in view of \eqref{LindX}  and
$\E \vep(\boldsymbol{0})^2 < \infty $. Thus, ${\cal S}_{\lambda, \gamma} \limd N(0, {\cal R}_\gamma) $, in  other words, the distribution of $ (d^{-1}_{\lambda,\gamma} S_{\lambda,\gamma}(\mbx_j);  j =1, \cdots, m ) $
tends to a Gaussian distribution on $\R^m $ with mean zero and covariance matrix  $(R_\gamma (\mbx_j, \mbx_{j'}))_{j,j=1, \cdots, m} $, proving
the proposition. \hfill $\Box$

\smallskip

Let $\upsilon_i >0, i=1,2, \Upsilon := \frac{1}{\upsilon_1} + \frac{1}{\upsilon_2}$ and $\rho(\mbx), \rho_p(\mbx) $ be as in \eqref{rhosp}, \eqref{rhosp2}.  We note the elementary inequality:
for any $p>0$  there exist constants $0 < C_1 \le C_2 < \infty $ such  that
\begin{equation}
C_-  \rho(\mbx) \le \rho_p(\mbx) \le C_+ \rho(\mbx), \qquad (\forall \, \mbx \in \R^2).
\end{equation}
We also use the  fact \cite{pils2017} that for any $\delta >0, w >0 $
\begin{equation}
\int_{|\mbx| <  \delta} \rho(\mbx)^{-w}  \d \mbx < \infty \ \Longleftrightarrow \ w < \Upsilon,
\qquad \int_{|\mbx| >  \delta} \rho(\mbx)^{-w}  \d \mbx < \infty \ \Longleftrightarrow \ w > \Upsilon.
\end{equation}
Following \cite{pils2020}, we say that a measurable function  $L : \R^2_0 \to \R$ is \emph{generalized invariant}
if the function $ \mbx \mapsto L(\lambda^{1/\upsilon_1}x_1, \lambda^{1/\upsilon_2} x_2)$ does not depend on $\lambda >0$.
Every generalized invariant function $L$ can be represented as
\begin{equation}\label{tildeL}
L(\mbx) = \tilde L (x_1/\rho(\mbx)^{1/\upsilon_1},
x_2/\rho(\mbx)^{1/\upsilon_2}), \quad \mbx \in \R^2_0,
\end{equation}
where $\tilde L$ is the restriction of $L$ to ${\cal S}_1 := \{\mbx \in \R^2_0: \rho(\mbx) = 1\}$.

\noi  {\bf Assumption (F)$_{\rm  LRD}
$} \ The spectral density $f$  satisfies
\begin{equation}\label{Fsp1}
f(\mbx)\  =  \  
\rho(\mbx)^{-1}L(\mbx) \big(1 + o(1)\big), \quad |\mbx| \to 0, 
\end{equation}
where $\rho(\mbx)$ as in \eqref{rhosp}, $\Upsilon > 1 $  
and $L(\mbx) = L(-\mbx), \mbx \in  \R^2_0$ is a strictly positive continuous generalized invariant function.
Moreover, $f$ is bounded  on $\{ \mbx \in \Pi^2: |\mbx| > \delta\}, $ for any $\delta >0$.

\smallskip

\noi  {\bf Assumption (F)$_{{\rm LRND},i}
$} \ $(i=1,2)$ \ The spectral density $f$  satisfies Assumption (F)$_{\rm  LRD}$ except that strict positivity of $L(\mbx) $ is replaced by 
\begin{equation} \label{Lmu}
\tilde L(\mbx) = \ell_i |x_i|^{\mu_i} (1 + o(1)) 
\qquad (x_i \to 0, \, \mbx \in {\cal S}_1)
\end{equation}
for some $\mu_i >0, \ell_i >0 $.

\medskip

\noi  {\bf Assumption (F)$_{\rm ND}$} \ The spectral density $f(\mbx), \mbx \in \Pi^2$ is a bounded continuous function such that
\begin{equation}\label{Fsp2}
f(\mbx)\  = \  \rho(\mbx) L(\mbx) 
\big(1 + o(1)\big), \quad |\mbx| \to  0, \quad \text{where}  \quad 
\end{equation}
where $\rho(\mbx)$ as in \eqref{rhosp}, 
and $L(\mbx), \mbx \in  \R^2_0$ is as in \eqref{Fsp1}.

With generalized homogeneous functions $ \rho(\mbx)^{-1}L(\mbx), \rho(\mbx)L(\mbx) $ in \eqref{Fsp1}, \eqref{Fsp2} we associate Gaussian RFs
$V_{0,1/\rho},  V_{0,\rho}$ by
\begin{eqnarray}  \label{V0G}
\hskip-.5cm &V_{0,1/\rho}(\mbx) := \int_{\R^2}\prod_{j=1}^2 \frac{1-\e^{{\i} u_j x_j}}{\i u_j}\, \sqrt{L(\mbu)/\rho(\mbu)}\, Z(\d \mbu),
\quad V_{0,\rho}(\mbx) := \int_{\R^2}\prod_{j=1}^2 \frac{1-\e^{{\i} u_j x_j}}{\i u_j}\, \sqrt{L(\mbu)\rho(\mbu)}\, Z(\d \mbu), 
\end{eqnarray}
$\mbx \in \R^2_+, $ where $\{Z(\d \mbu)\}$ is a complex-valued Gaussian noise with zero mean and $\E |Z(\d \mbu)|^2 = \d \mbu$.  The existence
of $ V_{0,1/\rho}, V_{0,\rho}$ will  be established in the following sections.

\begin{definition} Let $V = \{ V(\mbx); \mbx \in \R^2_+ \} $ be a RF and
$\gamma >0, H\ge 0, H_i \ge 0, i=1,2 $.  We say that

\noi (i) $V$ is   {\it $(\gamma, H)$-self-similar} (SS) if
\begin{equation} \label{Vss}
\{  V(\lambda x_1, \lambda^\gamma x_2); \mbx  \in \R^2_+\}
\eqfdd \{ \lambda^{H} V(\mbx); \mbx \in \R^2_+\}, \qquad \forall \lambda >0.
\end{equation}

\smallskip

\noi  (ii) $V$ is  {\it $(H_1, H_2)$-multi-self-similar} (MSS) is
\begin{equation} \label{Vmss}
\{  V(\lambda_1 x_1, \lambda_2 x_2); \mbx  \in \R^2_+\}
\eqfdd \{ \lambda_1^{H_1} \lambda_2^{H_2} V(\mbx); \mbx \in \R^2_+\}, \qquad \forall \lambda_1, \lambda_2 >0.
\end{equation}

\end{definition}

$(\gamma, H)$-SS  property is a particular case of the operator self-similarity  property introduced in \cite{bier2007} and corresponding
to scaling $\mbx \to \lambda^E \mbx $  with diagonal matrix $E = {\rm diag}(1, \gamma)$. Particularly,
$(1,H)$-SS property  coincides with the usual $H$-SS property for RFs on $\R^2_+$ \cite{samo2016}.  $(H_1, H_2)$-MSS property was introduced in
\cite{gent2007}. It implies $(\gamma, H)$-SS  property with  any $\gamma >0$ and $H = H_1 + \gamma H_2 $. Under mild additional assumptions,
scaling limits $V_\gamma^X$ in \eqref{VgammaD}  are $(\gamma, H)$-SS RFs and the normalization $d_{\lambda,\gamma} $ is regularly
varying at infinity with index $H$ \cite{ps2016}.

\begin{definition} 
\label{defFBS}
	(Standard) {\it Fractional Brownian Sheet} (FBS) $B_{H_1,H_2}= \{ B_{H_1,H_2} (\mbx); \mbx \in \R^2_+  \}$ with $(H_1,H_2)
	\in [0,1]^2$ 
is defined as a Gaussian process with zero-mean and covariance function $\E B_{H_1,H_2} (\mbx) B_{H_1,H_2} (\mby) =
\prod_{i=1}^2 r_{H_i} (x_i,y_i)$, $\mbx, \mby \in \R^2_+$, where for $x,y \in \R_+$,
	\begin{eqnarray} \label{FBScov}
	r_H (x,y)&:=&\frac{1}{2}\begin{cases}
x^{2H} + y^{2H} - |x-y|^{2H},  & 0<H\le 1, \\
2,  &H=0,  x = y,  \\
1, &H=0,  x \ne y.
\end{cases}
	\end{eqnarray}
\end{definition}

Usually, FBS is defined for $(H_1,H_2) \in (0,1]^2 $ or even  $(H_1,H_2) \in (0,1)^2 $, see \cite{aya2002}. Extension of FBS to  $H_1 \wedge H_2 = 0 $
was introduced in \cite{sur2020}.
\eqref{FBScov} implies that the restriction of  FBS  $B_{H_1,H_2}$ to horizontal/vertical line agrees with fractional Brownian motion (FBM) $B_H=\{B_H(x); x\in \R_+\}$ with covariance function $r_H(x,y)$ and
the corresponding Hurst parameter $H = H_i \in [0,1]$, $i=1,2 $. Note $r_0(x,y) = \lim_{H \downarrow 0} r_H(x,y)$.
FBM $B_0$ is  $0$-SS  and an extremely singular (non-measurable) process, see \cite[pp.256--257]{samo2016}.
It can be represented as $B_0 \eqfdd \{\frac{1}{\sqrt{2}}(W(x) - W(0)); x \in \R_+ \}$, where $W(x)$, $x \in [0,\infty)$, is (uncountable) family of \emph{independent} $N(0,1)$ r.v.s.
The above $B_0$ is  different from the 'regularized' FBM with $H=0$ defined in \cite[p.2985]{fyod2016}. 
FBS $B_{H_1, H_2}$ with $H_1\wedge H_2 =0$ and their $\alpha$-stable extensions appeared in limit theorems for RFs \cite{sur2020, pils2021}.
We also recall spectral representation of FBS with $(H_1,H_2) \in (0,1)^2 $:
\begin{eqnarray}\label{BHH}
B_{H_1,H_2}(\mbx)
&=&\kappa^{-1} \int_{\R^2} \prod_{j=1}^2 \frac{1- \e^{\i x_j u_j}}{\i u_j |u_j|^{H_j -\frac{1}{2}}} Z(\d \mbu)
\end{eqnarray}
where
\begin{eqnarray}\label{kappa}
\kappa^2&:=&\int_{\R^2} \prod_{j=1}^2 \frac{|1- \e^{\i u_j}|^2}{|u_j|^{2H_j +1}} \d \mbu \ = \ \prod_{j=1}^2 \frac{\pi}{H_j \Gamma(2H_j) \sin (H_j \pi)}
\end{eqnarray}
and $Z(\d \mbu) $ is the same Gaussian white noise as in \eqref{V0G}; see \cite{aya2002}. For $H_1 \vee H_2 =1 $, FBS is a line RF of the form
$B_{1,  H_2}(\mbx) = x_1 B_{H_2} (x_2) $ $(0< H_2 < 1)$,   $B_{H_1, 1}(\mbx) = x_2 B_{H_1} (x_1) $ $(0< H_1 < 1), B_{1,1}(\mbx) = x_1 x_2 Z $, where
$B_H $ is  FBM and $Z  \sim N(0,1)$. For $H_1 \wedge H_2 = 0$, FBS allow a construction of finite-dimensional  distributions
via independent FBM or uncountable family  of independent Gaussian variables \cite{sur2020}.

\section{Long-range dependence}

Define
\begin{eqnarray}\label{H+}
H^+_1&:=& \frac{1}{2}(1 + (\upsilon_1 \wedge 1)),  \quad   H^+_2  := \frac{1}{2}(1 + \upsilon_2 - \frac{\upsilon_2}{\upsilon_1 \vee 1}), \\
H^-_1&:=& \frac{1}{2}(1 + \upsilon_1 - \frac{\upsilon_1}{\upsilon_2 \vee 1}), \hskip.5cm
H^-_2:=\frac{1}{2}(1 + (\upsilon_2 \wedge 1)), \nn
\end{eqnarray}
Note $H^+_1 = 1 $ for $\upsilon_1 \ge 1 $, $H^+_2 = 1/2 $ for $\upsilon_1 \le 1 $. Analogous relations are satisfied by $H^-_i, i=1,2$.

\begin{theorem} \label{thmGaussLRD} Let $X$ in \eqref{Sgamma}
be a stationary linear RF on $\Z^2$ in \eqref{Xlin} with spectral density $f$ satisfying Assumption  (F)$_{LRD}$. Then:

\begin{itemize}

\item The scaling limits in \eqref{VgammaD} exist for any $\gamma >0$ and satisfy  \eqref{Vtrans}  with $\gamma_0 = \frac{\upsilon_1}{\upsilon_2} $ and
the unbalanced limits are given by
\begin{eqnarray}\label{VB}
V_+
&:=&\kappa_+
\begin{cases}
B_{H^+_1, 1/2}\, &\upsilon_1 \le 1, \\
B_{1, H^+_2},  
&\upsilon_1\ge 1,
\end{cases} \qquad
V_-
\ :=\  \kappa_-
\begin{cases}
B_{1/2, H^-_2} 
&\upsilon_2 \le 1, \\
B_{H^-_1, 1},
&\upsilon_2\ge 1,
\end{cases}
\end{eqnarray}
where $H^\pm_i,i=1,2 $ as in \eqref{H+}  and the
constants $\kappa_\pm >0 $ in \eqref{VB} are
defined in \eqref{kappa1}, \eqref{kappa2} (or can be determined from
these expressions by symmetry).

\item  The well-balanced limit $V_{0} = V_{0, 1/\rho}$ is given in \eqref{V0G}.

\item  The normalization in \eqref{VgammaD} is given by
$d_{\lambda, \gamma} := \lambda^{H(\gamma)} (\log_+ \lambda)^{1/2}$ if
$\gamma > \gamma_0 $ and $\upsilon_1   = 1$ or $\gamma < \gamma_0 $ and $\upsilon_2   = 1$, and
$d_{\lambda, \gamma} := \lambda^{H(\gamma)} $ in the remaining cases, where $H(\gamma)$ is the weighted linear combination with weights
$(1, \gamma)$ of the Hurst indices of FBM on the r.h.s. \eqref{VB}, viz.,
\begin{eqnarray} \label{Hgamma}
H(\gamma)= \begin{cases}H^+_1 + \gamma H^+_2, &\gamma \ge \gamma_0,  \\
H^-_1 + \gamma H^-_2, &\gamma \le \gamma_0.
\end{cases}
\end{eqnarray}

\item The RF $X$ exhibits scaling transition at $\gamma_0$.

\end{itemize}
\end{theorem}

\noi {\bf Proof.} The simple bound $\big|\sum_{\mbt \in K_{[\lambda x_1,\lambda^{\gamma} x_2]}} a(\mbt - \mbu)\big|
\le C \lambda^{(1+ \gamma)/2} \big( \sum_{\mbt \in  \Z^2} a(\mbt)^2\big)^{1/2} \le  C\lambda^{(1+\gamma)/2} $ reduces the criterion
\eqref{LindX} to $H(\gamma) >   (1+\gamma)/2$.  This follows from \eqref{Hgamma} since $H^\pm_i \ge 1/2, i=1,2 $.
Whence, for \eqref{VgammaD} it suffices to prove the convergence of covariance functions in   \eqref{Rlim}.
The latter relation 
is essentially proved in \cite{ps2015}
except for the cases $\upsilon_1 = 1  $ or $\upsilon_2 = 1 $. 
The subsequent proof is limited to these cases. 
Note that the definitions in \eqref{VB}, \eqref{Hgamma} given by different expressions agree for
$\upsilon_i =1, i=1,2 $ and/or $\gamma = \gamma_0$.   By symmetry, it suffices to discuss the case
$\upsilon_1 =1 $ only and consider the limits in \eqref{Rlim}
for $\gamma > \gamma_0, \gamma < \gamma_0, $ and  $\gamma = \gamma_0 $ separately.
W.l.g.  $\lambda > 1 $.

\smallskip

\noi {\it Case $\gamma > \gamma_0 = 1/\upsilon_2 \, (\upsilon_1 = 1)$.} In this case,
$V_+ = \kappa_+  B_{1,1/2 }$ and $d^2_{\lambda,\gamma} = \lambda^{2 + \gamma}  \log \lambda $,
according to  \eqref{VB}. In the integral in \eqref{Rcov}, change the variables:
$u_1 \to u_1/\lambda^{\gamma  \upsilon_2},  u_2 \to u_2/\lambda^{\gamma} $ and note that
\begin{eqnarray} \label{Dlim}
&&\lambda^{-1} D_{[\lambda x_1]}(u_1/\lambda^{\gamma  \upsilon_2}) \to x_1, \qquad
\lambda^{-\gamma} D_{[\lambda^\gamma x_2]}(u_2/\lambda^{\gamma}) \to \frac{1 -  \e^{\i x_2 u_2}}{- \i u_2}, \\
&&f_\lambda(\mbu) :=  \lambda^{-\gamma \upsilon_2} f(u_1/\lambda^{\gamma \upsilon_2}, u_2/\lambda^{\gamma}) \to L(\mbu)/\rho(\mbu)  =:
f_0(\mbu)\nn
\end{eqnarray}
point-wise for any $\mbu \in \R^2_0 $ as $\lambda \to \infty$,
due to $\gamma > \gamma_0$.  Moreover, \eqref{Fsp1}   implies
\begin{equation}\label{fgbdd}
|f_\lambda (\mbu) - f_0(\mbu)| \le \frac{1}{|u_1| + |u_2|^{\upsilon_2}}
\delta \big( \frac{u_1}{\lambda^{\gamma \upsilon_2}}, \frac{u_2}{\lambda^{\gamma}}\big),
\end{equation}
where $\delta(\mbu)$ is a bounded function tending to  0 as $\mbu \to \0$.
Thus,
$R_{\lambda, \gamma}(\mbx, \mby)= \sum_{i=1}^2 R^{(i)}_{\lambda, \gamma}(\mbx, \mby)$, 
where $R^{(i)}_{\lambda, \gamma}(\mbx, \mby) := \frac{1}{\log \lambda }  \int_{\lambda^{\gamma \upsilon_2} \Pi
\times \lambda^\gamma \Pi} G^{(i)}_\lambda (\mbu) \d \mbu, i=1,2$ with
\begin{eqnarray} \label{G1u}
G^{(1)}_\lambda (\mbu)
&:=&\frac{D_{[\lambda x_1]}(u_1/\lambda^{\gamma \upsilon_2})
\overline{D_{[\lambda y_1]}(u_1/\lambda^{\gamma \upsilon_2})}
}{\lambda^2} \cdot
\frac{D_{[\lambda^\gamma x_2]}(u_2/\lambda^{\gamma})
\overline{D_{[\lambda^\gamma y_2]}(u_2/\lambda^{\gamma})}}{\lambda^{2\gamma}}\cdot  f_0(\mbu) \\ 
&\to& x_1 y_1 \cdot \frac{(1 -  \e^{\i x_2 u_2}) (1 -  \e^{-\i y_2 u_2})}{|u_2|^2} \cdot f_0(\mbu) =: G(\mbu), \nn
\end{eqnarray}
and $G^{(2)}_\lambda (\mbu)$ analogously defined with $ f_0(\mbu) $ replaced by
$f_\lambda (\mbu) - f_0(\mbu)$.

Let us prove $R^{(2)}_{\lambda, \gamma}(\mbx, \mby) \to 0$. From elementary bound $|D_n(x)| \le C n/(1+ n|x|), x \in \Pi $
and \eqref{fgbdd} we get
\begin{eqnarray}\label{R2prime}
|R^{(2)}_{\lambda, \gamma}(\mbx, \mby)|
&\le&
 \frac{C}{\log \lambda }  \int_{\lambda^{\gamma \upsilon_2} \Pi
\times \lambda^\gamma \Pi}  \frac{ \d \mbu} {(1+ u_2^2)(|u_1| + |u_2|^{\upsilon_2})}
\delta \big( \frac{u_1}{\lambda^{\gamma \upsilon_2}}, \frac{u_2}{\lambda^{\gamma}} \big)\nn \\
&=& \frac{C}{\log \lambda }\int_{\R}  \frac{\d u_2}{1+ u_2^2} \int_{|u_1| \le \lambda^{\gamma \upsilon_2}/|u_2|^{\upsilon_2}}
\frac{\d u_1}{|u_1| + 1}\delta \big( \frac{u_1  |u_2|^{\upsilon_2}}{\lambda^{\gamma \upsilon_2}}, \frac{u_2}{\lambda^{\gamma}} \big).
\end{eqnarray}
Split the integration region on the r.h.s. of \eqref{R2prime} into three parts corresponding
to $\{\frac{|u_2|}{\lambda^\gamma} \ge \epsilon\}$, $
\{\frac{|u_2|}{\lambda^\gamma} <  \epsilon, |u_1| (\frac{|u_2|}{\lambda^\gamma})^{\upsilon_2} <  \epsilon\}$,
and $  \{\frac{|u_2|}{\lambda^\gamma} < \epsilon, |u_1| (\frac{|u_2|}{\lambda^\gamma})^{\upsilon_2} \ge \epsilon\}$,
$T_1, T_2 $ and $T_3$, say,
so that $ |R^{(2)}_{\lambda, \gamma}(\mbx, \mby)| \le C (\log \lambda)^{-1} \sum_{i=1}^3 \int_{T_i} \dots
$  $=: (\log \lambda)^{-1} \sum_{i=1}^3 R_i $.
The inner integral in \eqref{R2prime} on $T_1$ is bounded by $C \log (1/\epsilon)^{\upsilon_2}) $ implying $R_1 \to 0 \ (\lambda \to \infty, \, \forall \epsilon  >0)$.  Next,
by the property $\delta(\mbu) \to 0$, for any $\delta' >0 $  there exists $\epsilon >0$ s.t.
$R_2 \le \delta' (\log \lambda)^{-1} \int_{\R}   (1+ u^2_2)^{-1} \log (\lambda^\gamma/|u_2|)^{\upsilon_2} \d u_2  \le C \delta' $ vanishes
with $\epsilon \to 0 $ uniformly in $\lambda \ge 2$. Finally, for any $\epsilon >0$ fixed,
$R_3 \le C (\log \lambda)^{-1} \int_{\R}   (1+ u^2_2)^{-1} \big\{
\log (\lambda^\gamma/|u_2|)^{\upsilon_2} - \log \epsilon (\lambda^\gamma/|u_2|)^{\upsilon_2} \big\} = O(1/\log \lambda) \to 0$ as
 $\lambda \to \infty $, thus ending the proof of  $R^{(2)}_{\lambda, \gamma}(\mbx, \mby) \to 0$.

Next, we prove the limits
\begin{equation}\label{R1lim}
\lim_{\lambda \to \infty} R^{(1)}_{\lambda, \gamma}(\mbx, \mby) =  \lim_{\lambda \to \infty} \tilde R_{\lambda, \gamma}(\mbx, \mby) =
\kappa_+^2  \E [B_{1, 1/2}(\mbx) B_{1, 1/2}(\mby)] = \kappa_+^2 x_1 y_1 (x_2 \wedge y_2)
\end{equation}
where
\begin{eqnarray}
\tilde R_{\lambda, \gamma}(\mbx, \mby)&:=&\frac{1}{\log \lambda }  \int_{\lambda^{\gamma \upsilon_2} \Pi
\times \lambda^\gamma \Pi} G (\mbu) \d \mbu \nn \\
&=&x_1 y_1 \int_{|u_2| \le \lambda^\gamma \pi}  \frac{(1 -  \e^{\i x_2 u_2}) (1 -  \e^{-\i y_2 u_2})}{|u_2|^2} \d u_2 \times \frac{1}{\log \lambda }\int_{|u_1| \le \lambda^{\gamma \upsilon_2} \pi} f_0(\mbu) \d u_1. \label{R2lim}
\end{eqnarray}
Consider the second limit in \eqref{R1lim}.
By definition, $f_0 (\mbu) = L\big(\frac{u_1}{|u_1| + |u_2|^{\upsilon_2}},
\frac{u_2}{(|u_1| + |u_2|^{\upsilon_2})^{1/\upsilon_2}} \big) (|u_1| + |u_2|^{\upsilon_2})^{-1} $
behaves as $ L(1,0) |u_1|^{-1} $ when $|u_1| \to \infty $, indeed, for any $u_2 \ne 0$
\begin{equation}\label{R3lim}
|u_2|^{\upsilon_2} (|u_1| + 1) f_0(u_1 |u_2|^{\upsilon_2},u_2) =  L\Big( \frac{u_1 + 1}{|u_1| + 1},  \frac{{\rm sgn}(u_2)}{|u_1| + 1}\Big)
\ \to \ L(1,0),  \qquad |u_1| \to \infty
\end{equation}
by the conditions on $L(\cdot)$ in the theorem. Then, the last term on the r.h.s. of \eqref{R2lim} can be rewritten as the sum of two terms: the first term $\frac{L(1,0)}{\log \lambda }\int_{|u_1| \le (\lambda^\gamma/|u_2|)^{\upsilon_2} \pi} (|u_1|+1)^{-1}  \d u_1
\to 2 L(1,0) \gamma \upsilon_2 \ (\forall u_2 \ne 0)$, and the second term
\begin{equation}\label{R4lim}
\frac{1}{\log \lambda }\int_{|u_1| \le (\lambda^\gamma/|u_2|)^{\upsilon_2} \pi} \frac{\d u_1}{|u_1| + 1}
\Big(L\big( \frac{u_1 + 1}{|u_1| + 1},  \frac{{\rm sgn}(u_2)}{|u_1| + 1}\big) - L(1,0)\Big) \ \to \  0 \qquad (\forall \, u_2 \ne 0)
\end{equation}
in view of \eqref{R3lim} and boundedness of $L(\cdot)$. The proof of the second limit in \eqref{R1lim} follows from
\eqref{R2lim}-\eqref{R4lim} and the DCT using the  bound $(\log \lambda)^{-1} \big|\int_{|u_1| \le \lambda^{\gamma \upsilon_2} \pi} g_+(\mbu) \d u_1 \big|
\le C (\log \lambda)^{-1} \int_{|u_1| \le \lambda^{\gamma \upsilon_2} \pi} (|u_1| + |u_2|^{\upsilon_2})^{-1} \d u_1 \le C(1 +  |\log (|u_2|)|)    $
and the integrability of $ (1 + |\log (|u_2|)|)(1+ u_2^2)^{-1}  $ on $\R$.

It remains to prove that $R^{(1)}_{\lambda, \gamma}(\mbx, \mby) - \tilde R_{\lambda, \gamma}(\mbx, \mby) \to 0$. W.l.g.,
take $\mbx = \mby = \1 $ and let $R^{(1)}_{\lambda, \gamma}(\1, \1) - \tilde R_{\lambda, \gamma}(\1, \1) =:  R'_{\lambda, \gamma}$. 
Then
\begin{eqnarray}\label{R5lim}
|R'_{\lambda,\gamma}|
&\le&
\frac{1}{\log \lambda} \int_{\lambda^{\gamma \upsilon_2} \Pi \times \lambda^\gamma \Pi}
\Big| \frac{|D_{[\lambda]}(u_1/\lambda^{\gamma  \upsilon_2})D_{[\lambda^\gamma]}(u_2/\lambda^{\gamma})|^2}{\lambda^{2 + 2 \gamma}}
- \frac{|1 -  \e^{\i u_2}|^2}{u_2^2} \Big|  \frac{\d \mbu}{|u_1| + |u_2|^{\upsilon_2}}\nn \\
&=&\frac{1}{\log \lambda} \int_{\lambda^{\gamma \upsilon_2} \Pi \times \lambda^\gamma \Pi}
\frac{\delta_\lambda (\mbu) \d \mbu }{(1 + u_2^2) (|u_1| + |u_2|^{\upsilon_2})} \ = \ o(1)
\end{eqnarray}
since $\delta_\lambda (\mbu) \to  0$ uniformly in $\mbu $ in the last integral,
see  \eqref{Dlim}, and the r.h.s. of \eqref{R5lim} with $\delta_\lambda (\mbu) $ replaced
by 1 is bounded.
This ends the proof of  \eqref{R1lim} with
\begin{equation}\label{kappa1}
\kappa_+^2
:=  2 L(1,0) \gamma \upsilon_2 \int_{\R}  |(1 -  \e^{\i u})/u|^2 \d u \ =  \ 4 \pi L(1,0) \gamma \upsilon_2
\end{equation}
and completes
the proof of  \eqref{Rlim} when $\gamma > \gamma_0$.

\smallskip

\noi {\it Case $\gamma < \gamma_0 = 1/\upsilon_2 \, (\upsilon_1 = 1)$.} There are three subcases:
(a) $\upsilon_2 < 1,$ (b) $\upsilon_2 > 1 $, and (c) $\upsilon_2 = 1 $.
Accordingly, $V_- = \kappa_-  B_{1/2, H^-_2}$  in subcase (a),
$V_- = \kappa_-  B_{H^-_1, 1} $ in subcase (b), and
$V_- = \kappa_-  B_{1/2, 1} $ 
in subcase (c).
Subcases  (a)  and (b) do not involve logarithmic normalization and are essentially treated in \cite[Thm.3.1]{ps2015}. Subcase (c)
is symmetric to the above case $\gamma >\gamma^X_0, \upsilon_1 = 1 $ with $(x_1, \upsilon_1) $ and $(x_2, \upsilon_2) $ exchanged, by
noting that the latter discussion applies  to any $\upsilon_2 >0$ including $\upsilon_2 = 1$.

\smallskip

\noi {\it Case $\gamma = \gamma_0 = 1/\upsilon_2 \, (\upsilon_1 = 1)$.} In this case,
$V_0 = V_{0, 1/\rho}$ is  given in \eqref{V0G} and the result follows from \cite{ps2015} with small changes, thus
ending the proof of Theorem \ref{thmGaussLRD}. \hfill $\Box$

\smallskip

\begin{remark} {\rm The squared asymptotic constant in \eqref{VB} for $\upsilon_1 \ne 1 $ takes the form
\begin{eqnarray} \label{kappa2}
\kappa^2_+
&=&\begin{cases}
L(1,0)\int_{\R^2} \prod_{j=1}^2 \big|(1-\e^{\i u_j})/u_j\big|^2 |u_1|^{-\upsilon_1} \d \mbu,  &\upsilon_1 < 1, \\
\int_{\R^2} \big|(1-\e^{\i u_2})/u_2\big|^2 L(\mbu)\rho(\mbu)^{-1}  \d \mbu,  &\upsilon_1 > 1,
\end{cases}
\end{eqnarray}
($\kappa^2_- $ is defined symmetrically with $u_1, u_2, \upsilon_1, \upsilon_2 $ exchanged by $u_2, u_1, \upsilon_2, \upsilon_1 $).
For $\upsilon_1 < 1$ the integral in  \eqref{kappa2} can be explicitly evaluated as
$\kappa^2_+ = L(1,0) (2\pi)^2/\Gamma (2+ \upsilon_1) \sin
((1+ \upsilon_1)\pi/2),
 $ see \eqref{kappa}.
}
\end{remark}

\section{Negative dependence}

This sec. describes anisotropic scaling  limits in \eqref{VgammaD} of linear
RFs satisfying Assumption (F)$_{\rm ND}$. Define
\begin{eqnarray}\label{H-}
H^+_1 := \frac{1}{2}(1 - (\upsilon_1 \wedge 1)), \qquad
H^-_2:= \frac{1}{2}(1 - (\upsilon_2 \wedge 1)), \qquad \gamma_0 := \frac{\upsilon_1 \wedge 1}{\upsilon_2 \wedge 1}. 
\end{eqnarray}
Note $H^+_1, H^-_2 \in [0,1/2)$ and $H^+_1 = 0$ (respectively,  $H^-_2 = 0$)
is equivalent to $\upsilon_1 \ge 1 $ (respectively,  $\upsilon_2 \ge 1$). We  also set
$H^+_2 = H^-_1 := 1/2 $.

\begin{theorem} \label{thmGaussND} Let $X$ in \eqref{Sgamma}
be a stationary linear RF on $\Z^2$ in \eqref{Xlin} with spectral density $f$
satisfying Assumption  (F)$_{ ND}
$.   Then:

\begin{itemize}

\item The scaling limits in \eqref{VgammaD} exist for any $\gamma >0$ and satisfy  \eqref{Vtrans}  with
$\gamma_0, H^\pm_i$ in \eqref{H-} and the unbalanced limits given by
\begin{eqnarray}\label{VBneg}
V_+
&:=&\kappa_+ B_{H^+_1, 1/2}, \qquad
V_-
\ :=\  \kappa_- B_{1/2, H^-_2}. 
\end{eqnarray}
The asymptotic constants $\kappa_\pm >0 $ 
are written in \eqref{kappa3}, \eqref{kappa4}, \eqref{kappa44}, \eqref{kappa5}, and \eqref{kappa6} (or can be determined from
these expressions by symmetry).

\item  The well-balanced limit is given by 
\begin{eqnarray}\label{V0neg}
V_0
&:=&\begin{cases}
V_{0, \rho}, &H^+_1 \wedge H^-_2 >0, \\
\kappa_+ B_{H^+_1,1/2} +  \kappa_- B_{1/2, H^-_2},  &H^+_1 \wedge H^-_2= 0,
\end{cases}
\end{eqnarray}
where $B_{H^+_1,1/2}$ and $B_{1/2, H^-_2}$ are independent and $V_{0,\rho}$ is defined in \eqref{V0G}.

\item  The normalization in \eqref{VgammaD} is given by
$d_{\lambda, \gamma} :=  \lambda^{H(\gamma)} (\log_+ \lambda)^{1/2} $ in the  cases
$\gamma \ge \gamma_0,  H^+_1 = 0$
or $\gamma \le \gamma_0,  H^-_2 = 0 $, and $d_{\lambda, \gamma} :=  \lambda^{H(\gamma)}  $ in the  remaining cases,
with
\begin{eqnarray} \label{HgammaN}
H(\gamma)&:=&
\begin{cases}H^+_1 + (\gamma/2), &\gamma \ge \gamma_0,  \\
(1/2) + \gamma H^-_2, &\gamma \le \gamma_0.
\end{cases}
\end{eqnarray}

\item The RF $X$ exhibits  scaling transition at $\gamma_0 $.

\end{itemize}
\end{theorem}

\noi {\bf Proof.} Similarly as in the proof of Theorem \ref{thmGaussLRD} we check the asymptotic gaussianity  criterion in
\eqref{LindX},
reducing the proof to the convergence in \eqref{Rlim} of covariance functions. Relation
$f(\mbu) = (2\pi)^2 |\widehat a(\mbu)|^2 $ and the assumptions on $f$ imply that $|\widehat a(\mbu)| \le C $ is  bounded. Whence,
with  $(\gamma_1, \gamma_2) := (1, \gamma)$, we get that
$\sup_{\mbs \in \Z^2}  \big|\sum_{\mbt \in K_{[\lambda x_1,\lambda^{\gamma} x_2]}} a(\mbt - \mbs)\big|
\le \int_{\Pi^2} \prod_{i=1}^2 \big|D_{[\lambda^{\gamma_i} x_i]} (u_i)\big| \, |\widehat a(\mbu)| \d \mbu
\le C  \int_{\Pi^2} \prod_{i=1}^2 \big|D_{[\lambda^{\gamma_i} x_i]} (u_i)\big| \, \d \mbu \le C (\log \lambda)^2 $ so that
\eqref{LindX} holds since $H(\gamma)$ in \eqref{HgammaN} satisfy $H(\gamma)  >0$.

The subsequent proof of \eqref{Rlim} is
split into parts I and II according to whether $\upsilon_i \neq 1, i=1,2, $ or
 $\upsilon_i =1 \, (\exists i=1,2)$.  In turn, each part is split into several cases and subcases 
which are treated separately.
Part II involves the logarithmic factor in the normalization and is more delicate.

\smallskip

\noi {\bf Part I.} \underline{\it Case $\upsilon_1 \vee \upsilon_2  < 1, \gamma_0=
\upsilon_1/\upsilon_2$. }  Then
$H^+_1 = (1 - \upsilon_1)/2 \in (0, 1/2),  H^-_2 = (1 - \upsilon_2)/2 \in (0,1/2)$. By change of variables:
$u_1 \to u_1/\lambda, u_2 \to u_2/\lambda^\gamma $ we rewrite \eqref{Rcov} as
$R_{\lambda, \gamma}(\mbx, \mby) = \int_{\lambda \Pi \times \lambda^\gamma \Pi} G_\lambda (\mbu)  \d \mbu $ where
\begin{eqnarray}  \label{Guv}
G_\lambda(\mbu)&:=&
\frac{D_{[\lambda x_1]}(u_1/\lambda) \overline{D_{[\lambda y_1]}(u_1/\lambda)}} {\lambda^2} \cdot
\frac{D_{[\lambda^\gamma x_2]}(u_2/\lambda^\gamma) \overline{D_{[\lambda^\gamma y_2]}(u_2/\lambda^\gamma)}} {\lambda^{2\gamma}} \cdot
f_\lambda (\mbu), \\
f_\lambda (\mbu)&:=&  f(u_1/\lambda, u_2/\lambda^\gamma)/\lambda^{2H(\gamma) - 1 - \gamma}. \nn
\end{eqnarray}
By \eqref{Fsp1} and the definition of $H(\gamma)$ in \eqref{HgammaN},
\begin{eqnarray} \label{Dlim}
&&\lambda^{-1} D_{[\lambda x_1]}(u_1/\lambda) \to \frac{1 -  \e^{\i x_1 u_1}}{- \i u_1}, \qquad
\lambda^{-\gamma} D_{[\lambda^\gamma x_2]}(u_2/\lambda^{\gamma}) \to \frac{1 -  \e^{\i x_2 u_2}}{- \i u_2}, \\
&&f_\lambda(\mbu) \to g_\gamma(\mbu) \ := \  \begin{cases}
f_0(\mbu), &\gamma = \gamma_0, \\
f_0(u_1,0), &\gamma > \gamma_0, \\
f_0(0,u_2), &\gamma < \gamma_0,
\end{cases}\nn
\end{eqnarray}
point-wise for any $\mbu \in \R^2_0$,
where $f_0 (\mbu) :=  L(\mbu) \rho(\mbu)     $, see \eqref{Fsp2}. Particularly, $f_0(u_1,0) = |u_1|^{-\upsilon_1}
L(1,0) $, $f_0(0,u_2) = |u_2|^{-\upsilon_2}
L(0,1), \mbu = (u_1,u_2) \in \R^2_0 $
 since $L(1,0) = L(-1,0), $  $ L(0,1) = L(0,-1)$.
 We also see from
  \eqref{BHH}, \eqref{V0G} that the covariance function of the limit RF with appropriately chosen $\kappa_\pm $
 writes as
  $\E V_\gamma (\mbx) V_\gamma (\mby) = \int_{\R^2} G(\mbu) \d \mbu$, where
\begin{equation}\label{Gdef}
G(\mbu) \ := \ \prod_{j=1}^2  \big(\frac{1- \e^{\i x_j u_j}}{\i u_j}\big)
\big(\frac{1- \e^{-\i y_j u_j}}{-\i u_j}\big) 
g_\gamma(\mbu)
\end{equation}
is the product of the corresponding limit functions in \eqref{Dlim} and $ G_\lambda (\mbu)\to G(\mbu), \mbu \in \R^2_0 $  according to
\eqref{Dlim}. The proof of $\int_{\R^2} G_\lambda (\mbu)  \d \mbu \to \int_{\R^2} G(\mbu) \d \mbu$ in all three cases
$\gamma > \gamma_0, \gamma < \gamma_0$, and $\gamma = \gamma_0 $ now follows from the dominating convergence theorem (DCT) using the bound
\begin{equation*}
|\lambda^{-1}D_{[\lambda x]}(\mbox{$\frac{u}{\lambda}$})| \le C x (1+ |([\lambda x]/\lambda ) u|) \le C/(1+ |u|),\qquad |u|< \lambda \pi,
 \end{equation*}
for any fixed $x \in \R, x \ne 0$, implying
\begin{eqnarray}\label{Gdom}
&|G_\lambda (\mbu)| \ \le \ C \prod_{i=1}^2 (1+ u_i^2)^{-1} \times \begin{cases}
\rho(\mbu), &\gamma = \gamma_0, \\
\rho(u_1,0),  &\gamma > \gamma_0, \\
\rho(0,u_2),  &\gamma < \gamma_0.
\end{cases}
\end{eqnarray}
Hence, the r.h.s. of \eqref{Gdom} denoted by $\bar G(\mbu)$ is integrable:
$\int_{\R^2} \bar G(\mbu) \d \mbu < \infty $.
The asymptotic constants $\kappa_\pm $
in \eqref{VBneg} are determined by $\int_{\R^2} |G(\mbu)|^2 \d \mbu = 1/\kappa^2_+ (\gamma > \gamma_0), = 1/\kappa^2_- (\gamma < \gamma_0)
$ for $\mbx = \mby = \1 $ in \eqref{Gdef} and
take a similar form as in \eqref{kappa2}; particularly,
\begin{eqnarray}\label{kappa3}
\kappa_+^2&:=&
L(1,0) \int_{\R^2}  |u_1|^{\upsilon_1} \prod_{j=1}^2 |(1-\e^{{\i} u_j})/u_j|^2  \d \mbu.
\end{eqnarray}

\smallskip

\noi \underline{\it Case $\upsilon_1 < 1 < \upsilon_2, \gamma_0=
\upsilon_1$. }

\smallskip

\noi {\it Subcase $\gamma > \upsilon_1$.} We  prove \eqref{Rlim} with the same   $
V_+ = \kappa_+ B_{H^+_1, 1/2}$ as in the previous case using
a modified argument, as follows. 
Split $f(\mbu) = f_1(\mbu) + \tilde f_1(\mbu), $ where
$f_1(\mbu) = f(u_1,0),   \tilde f_1(\mbu) := f(\mbu) - f(u_1,0)$ and accordingly,
$R_{\lambda, \gamma}(\mbx, \mby) = R_1(\mbx, \mby) + \tilde R_1(\mbx, \mby)$.
The convergence $R_1(\mbx, \mby) \to \kappa^2_+ \E B_{H^+_1,1/2}(\mbx) B_{H^+_1,1/2}(\mby) $ for any $\gamma >0$ follows
as in the case $\upsilon_1 \vee \upsilon_2  < 1$ above.  Whence, it suffices to prove $ \tilde R_1(\mbx, \mbx) \to 0$, or
\begin{eqnarray} \label{Jsub1}
J_\lambda := \int_{\Pi^2} |\tilde f_1(\mbu)|\, |D_{[\lambda x_1]}(u_1) D_{[\lambda^{\gamma} x_2]}(u_2)|^2
\d \mbu
&=&o(\lambda^{1- \upsilon_1 + \gamma}).  
\end{eqnarray}
We have $|\tilde f_1(\mbu)| \le C \sum_{k=1}^3 |h_{1,k}(\mbu)|  $, where
$h_{1,1} (\mbu) := \rho(\mbu) -  \rho(u_1,0) = |u_2|^{\upsilon_2},
h_{1,2}(\mbu) := \rho(u_1,0)(L(\mbu) - L(u_1,0)) $  and $h_{1,3}(\mbu) = \rho(\mbu) \delta (\mbu) $, where $\delta (\mbu) $ is a bounded function
tending  to  $0$ as $|\mbu| \to 0 $. Accordingly, $J_\lambda \le C \sum_{k=1}^3 J_{\lambda, k} $ where $J_{\lambda, k}  $ is defined
as in \eqref{Jsub1} with $|\tilde f_1(\mbu)| $ replaced by $|h_{1,k}(\mbu)|$. Here,
\begin{eqnarray}\label{Jsub11}
J_{\lambda,1}
&\le&\int_{\Pi^2}|u_2|^{\upsilon_2}   \frac{|1 - \e^{\i [\lambda x_1] u_1}|^2}{ |u_1|^2 |u_2|^2}  \d \mbu \
\le \lambda \int_{\R \times \Pi}  \frac{|1 - \e^{\i u_1 x_1}|^2}{ |u_1|^2} |u_2|^{\upsilon_2 -2} \d \mbu  \le  C \lambda
= o(\lambda^{1- \upsilon_1 + \gamma})
\end{eqnarray}
since $\gamma > \upsilon_1 $. Next, 
\begin{eqnarray}
J_{\lambda,2} \label{Jsub12}
&\le&\int_{\Pi} |u_1|^{\upsilon_1} \frac{|1 - \e^{\i [\lambda x_1] u_1 }|^2}{ |u_1|^2} \d u_1 \times
\int_{\Pi} |L(u_1,u_2)- L(u_1,0)|   \frac{|1 - \e^{\i [\lambda^\gamma x_2] u_2}|^2}{|u_2|^2}  \d u_2 \\
&\le&C\lambda^{1-\upsilon_1 + \gamma}
\int_{\R} |u_1|^{\upsilon_1} \frac{|1 - \e^{\i u_1 x_1}|^2}{ |u_1|^2} \d u_1  \times
\int_{\R} \big|L\big(\frac{u_1}{\lambda}, \frac{u_2}{\lambda^\gamma}\big) - L\big(\frac{u_1}{\lambda}, 0\big)\big|
 \frac{|1 - \e^{\i u_2 x_2}|^2}{ |u_2|^2} \d u_2 \nn \\
&=&o(\lambda^{1- \upsilon_1 + \gamma}) \nn
\end{eqnarray}
by the DCT and the fact that $\big|L\big(\frac{u_1}{\lambda}, \frac{u_2}{\lambda^\gamma}\big) - L\big(\frac{u_1}{\lambda}, 0\big)\big| $
is bounded and tends to  0 for any  $\mbu \in \R^2_0 $  by the continuity of $\tilde L$ in  \eqref{tildeL} and the fact that
$\gamma > \upsilon_1 >  \upsilon_1/\upsilon_2 $.  A similar argument also entails
$J_{\lambda,3} = o(\lambda^{1- \upsilon_1 + \gamma})$,
proving  \eqref{Jsub1} and the convergence \eqref{VgammaD} in the subcase $\gamma > \upsilon_1$.

\smallskip

\noi {\it Subcase $\gamma < \upsilon_1$.} We prove  \eqref{Rlim} with $V_- = \kappa_-  B_{1/2, 0}$,
$d^2_{\lambda,\gamma}  = \lambda $ and $\kappa_-$ in \eqref{kappa4}.
Split $f(\mbu) = f_2(\mbu) + \tilde f_2(\mbu), \tilde f_2(\mbu) := f(\mbu) - f(0, u_2)$ and accordingly,
$R_{\lambda, \gamma}(\mbx, \mby) =  R_2(\mbx, \mby) + \tilde R_2(\mbx, \mby)$. It suffices to prove
\begin{eqnarray}\label{Jsub2}
R_2(\mbx, \mby) \to \kappa_-^2 (x_1 \wedge y_1)\times \frac{1}{2}(1 +   \I(x_2 = y_2)) \quad \text{and}
\quad \tilde R_2(\mbx, \mbx) \to  0.
\end{eqnarray}
To show the first relation in \eqref{Jsub2}, note that
\begin{eqnarray*}
R_2(\mbx, \mby)
&=&\frac{1}{\lambda}
 \int_{\Pi} D_{[\lambda x_1]}(u) \overline{D_{[\lambda y_1]}(u)} \d u \times
\int_{\Pi} f(0,v)    D_{[\lambda^{\gamma} x_2]}(v) \overline{D_{[\lambda^{\gamma} y_2]}(v)} \d v =: J_1 \times J_2,
\end{eqnarray*}
where
\begin{eqnarray*}
J_1&\to&
\int_{\R} \frac{(1- \e^{\i x_1 u})(1-  \e^{-\i y_1 u})}{|u|^2}  \d u
= (x_1 \wedge y_1) \int_{\R} \frac{|1-\e^{\i u}|^2}{|u|^2}  \d u
\end{eqnarray*}
follows by the DCT, and
\begin{eqnarray*}
J_2&=&
\int_{\Pi}\frac{f(0,v)}{|1- \e^{\i v}|^2} \big(1 - \e^{\i [\lambda^\gamma x_2] v}
-  \e^{-  \i [\lambda^\gamma y_2] v} +  \e^{\i ([\lambda^\gamma x_2] - [\lambda^\gamma y_2]) v}\big) \d v \\
&\to&\int_{\Pi}\frac{f(0,v)}{|1- \e^{\i v}|^2} \d v \times
\begin{cases}2, &x_2 = y_2, \\
1, &x_2 \ne y_2,
\end{cases}
\end{eqnarray*}
by the Lebesgue-Riemann theorem \cite[Thm.1.2]{ste1971}
and the integrability of $f(0,v) |1- \e^{\i v}|^{-2} $.
Whence, the asymptotic constant in \eqref{Jsub2} equals
\begin{equation} \label{kappa4}
\kappa_-^2
\ : =\  2 \int_{\R} |(1- \e^{\i u})/u|^{2}  \d u  \times \int_{\Pi} f(0,v) |1- \e^{\i v}|^{-2} \d v.
\end{equation}
To show the second relation in \eqref{Jsub2}, similarly as in the  subcase $\gamma > \upsilon_1$ above
write   $|\tilde f_2(\mbu)| \le C \sum_{k=1}^3 |h_{2,k}(\mbu)|  $ and $\tilde R_2(\mbx, \mbx) \le C \sum_{k=1}^3 \tilde R_{3, k}(\mbx, \mbx)$ accordingly, where $h_{2,1} (\mbu) := \rho(\mbu) -  \rho(0,u_2) = |u_1|^{\upsilon_1}, $  $
h_{2,2}(\mbu) :=  $   $\rho(0,u_2)(L(\mbu) - L(0,u_2)) $  and $h_{2,3}(\mbu) := \rho(\mbu) \delta (\mbu) $, where $\delta (\mbu) $ is a bounded function
tending  to  $0$ as $|\mbu| \to 0 $.  Then
$\tilde R_{2,3}(\mbx, \mbx) = o(1)$, 
\begin{eqnarray}
|\tilde R_{2,1}(\mbx, \mbx)|
&\le&C\lambda^{-1}
\int_{\Pi^2} |u_1|^{\upsilon_1} \frac{|1- \e^{\i  [\lambda x_1] u_1}|^2 |1- \e^{\i  [\lambda^\gamma x_2] u_2}|^2}
{ |1- \e^{\i  u_1}|^2  |1- \e^{\i  u_2}|^2} \d \mbu \nn  \\
&\le&C \lambda^{\gamma - \upsilon_1} \int_{\R} |u_1|^{\upsilon_1 -2} |1- \e^{\i  x_1 u_1}|^2 \d u_1 \times
\int_{\R}  |u_2|^{-2} |1- \e^{\i x_2 u_2}|^2 \d u_2 \nn \\
&=&O(\lambda^{\gamma - \upsilon_1}) = o(1) \label{Jsub3}
\end{eqnarray}
and
\begin{eqnarray}
|\tilde R_{2,2}(\mbx, \mbx)| 
&\le&C \int_{\Pi}  |u_2|^{\upsilon_2-2} \d u_2    \int_{\R} \big|L\big(\frac{u_1}{\lambda}, u_2) -
L(0, u_2)\big| |(1- \e^{\i x_1 u_1})/u_1|^2 \d u_1, 
= o(1), \label{Jsub4}
\end{eqnarray}
proving \eqref{Jsub2}. 

\smallskip

\noi {\it Subcase $\gamma  = \upsilon_1$.} We prove  \eqref{Rlim} with $V_0   = \kappa_+ B_{H^+_1,1/2}
+ \kappa_- B_{1/2,0} $ a sum of independent limits in subcases $\gamma > \upsilon_1 $ and $\gamma < \upsilon_1$.
Split $f(\mbu) = f_1(\mbu) +  f_2(\mbu) + \tilde f_{12}(\mbu),
\tilde f_{12}(\mbu) := f(\mbu) - f(u_1,0) - f(0, u_2)$ and, accordingly,
$R_{\lambda, \gamma}(\mbx, \mby) = R_1(\mbx, \mby)+ R_2(\mbx, \mby) + \tilde R_{12}(\mbx, \mby)$. Note $H^+_1 + (\upsilon_1/2) = 1/2 $.
The convergences
$R_1(\mbx, \mby) \to \kappa^2_+ \E B_{H^+_1,1/2} (\mbx, \mby) B_{H^+_1,1/2}(\mbx, \mby) $ and
$R_2(\mbx, \mby) \to \kappa^2_- \E B_{1/2,0} (\mbx, \mby) B_{1/2,0}(\mbx, \mby) $ were proved
in subcases  $\gamma > \gamma_0, \upsilon_i < 1, i=1,2 $ and $\gamma < \gamma_0,  \upsilon_1 < 1 < \upsilon_2 $, respectively. It remains to
show
\begin{eqnarray} \label{tiR12}
\tilde R_{12}(\mbx, \mbx)  \to 0.
\end{eqnarray}
Write $f(\mbu) = f_0(\mbu) T(\mbu)$, where $f_0(\mbu) = \rho(\mbu)  L(\mbu)$ and
$T(\mbu) :=  \frac{f(\mbu)}{f_0 (\mbu)}, \mbu \in \Pi^2 $ are a continuous functions,
see Assumption  (F)$_{\rm ND}$.
Then, $\tilde f_{12}(\mbu) = \sum_{i=1}^2 \tilde f_{12,i}(\mbu)$, where
\begin{eqnarray*}
\tilde f_{12,1}(\mbu) := |u_1|^{\upsilon_1} (\widetilde T(\mbu) - \widetilde T(u_1,0)), \qquad
\tilde f_{12,2}(\mbu) := |u_2|^{\upsilon_2} (\widetilde T(\mbu) - \widetilde T(0,u_2)), \quad \widetilde T(\mbu) := L(\mbu) T(\mbu).
\end{eqnarray*}
Accordingly,
$\tilde R_{12}(\mbx, \mbx) = \sum_{i=1}^2 \tilde R_{12,i}(\mbx, \mbx)$. Then
\begin{eqnarray}
\tilde R_{12,1}(\mbx,\mbx) 
&\le&\int_{\R^2} |u_1|^{\upsilon_1} \delta_{\lambda,1} (\mbu) 
\prod_{j=1}^2
|(1- \e^{\i u_j  x_j})/u_j|^2 \d \mbu,  \label{T1}
\end{eqnarray}
where $\delta_{\lambda,1} (\mbu) :=  \big|\widetilde T\big(\frac{u_1}{\lambda}, \frac{u_2}{\lambda^{\upsilon_1}} \big) - \widetilde T\big(\frac{u_1}{\lambda}, 0 \big)\big|\I(\mbu \in \lambda \Pi \times \lambda^{\upsilon_1} \Pi) $ is bounded  uniformly in $\lambda > 1 $ and
the integral on the r.h.s. with $\delta_{\lambda,1} (\mbu)$ replaced by 1 converges due to $\upsilon_1 < 1 $. Also note that
$\delta_{\lambda,1} (\mbu) \to 0$ as $T(\mbu), \mbu \in \Pi^2 $ is continuous and $L\big(\frac{u_1}{\lambda}, \frac{u_2}{\lambda^{\upsilon_1}} \big) - L\big(\frac{u_1}{\lambda}, 0 \big) \to 0$ due to $\upsilon_2 > 1 $. This proves \eqref{tiR12} for
$\tilde R_{12,1}(\mbx, \mbx)$ instead of $\tilde R_{12}(\mbx, \mbx)$. The proof of $\tilde R_{12,2}(\mbx, \mbx) \to 0$ is similar to
\eqref{Jsub4}, with $\big|L\big(\frac{u_1}{\lambda}, u_2) -
L(0, u_2)\big|$ there replaced by $\delta_{\lambda,2} (\mbu) := \big|\widetilde T\big(\frac{u_1}{\lambda}, u_2) -
\widetilde T (0, u_2)\big| \to 0 $  in view of the above mentioned properties of
$L(\mbu)$ and $T(\mbu)$.

\smallskip

\noi \underline{\it Case $\upsilon_i > 1, i=1,2,   \gamma_0= 1$. }  By symmetry, it suffices to  consider the  case  $\gamma \le 1 $.
Split $f(\mbu) =  f_1 (\mbu) +  f_2(\mbu) + \tilde f_{12}(\mbu), $  $R_{\lambda, \gamma}(\mbx, \mby) = R_1(\mbx, \mby)+ R_2(\mbx, \mby) + \tilde R_{12}(\mbx, \mby)$ 
as in the previous case. 
Then $R_2(\mbx, \mby) \to \kappa^2_- \E B_{1/2,0}(\mbx) B_{1/2,0}(\mby) $ as in
\eqref{Jsub2} while 
$R_1(\mbx, \mbx) \le C \lambda^{\gamma -1} \int_{\Pi} |u_1|^{\upsilon_1 -2} \d u_1 $  $
\times \int_{\R} |(1 - \e^{\i u_2 x_2})/u_2|^2 \d u_2 =
O(\lambda^{\gamma-1}) = o(1)  $ is  negligible when $\gamma < 1$; for  $\gamma = 1 $ we have
$R_1(\mbx, \mby) \to \kappa^2_+ \E B_{0,1/2}(\mbx) B_{0, 1/2}(\mby) $ analogously  to \eqref{Jsub2}  with
\begin{equation} \label{kappa44}
\kappa_+^2
\ : =\ 2\int_{\R} |(1- \e^{\i u})/u|^{2}  \d u  \times \int_{\Pi} f(v,0) |1- \e^{\i v}|^{-2} \d v.
\end{equation}
The proof of $\tilde R_{12}(\mbx, \mbx) \to 0$ for $\gamma \le 1$ follows  as in \eqref{Jsub4}.  This  ends
the proof of Part I.

\medskip

\noi {\bf Part II.} 
\underline{\it Case $\upsilon_1 =1, \upsilon_2 <1, \gamma_0 = 1/\upsilon_2  $. } 
Let first $\gamma > 1/\upsilon_2, d^2_{\lambda,\gamma} = \lambda^{\gamma} \log_+ \lambda $.
Split $f(\mbu) =  f_1(\mbu) + \tilde f_1(\mbu) $  
and
$R_{\lambda, \gamma}(\mbx, \mby) = R_1(\mbx, \mby) + \tilde R_1(\mbx, \mby) $ as in the case $\upsilon_1 < 1 < \upsilon_2$ above. Then
\begin{eqnarray}
R_1(\mbx, \mby) 
&\sim&\frac{1}{\log \lambda} \int_{\lambda \Pi} \frac{(1- \e^{\i x_1 u_1})(1- \e^{-\i y_1 u_1})}{|u_1|^2} \lambda f(u_1/\lambda, 0) \d u_1 \nn \\
&& \times \int_{\lambda^\gamma \Pi} \frac{(1- \e^{\i x_2 u_2})(1- \e^{-\i y_2 u_2})}{|u_2|^2} \d u_2 \  =: \  J_1 \times J_2,
\label{Jsub6}
\end{eqnarray}
where $J_2 \to  \kappa^2 (x_2 \wedge y_2),  \kappa^2 := \int_{\R} |(1-\e^{\i u})/u|^2 \d u $ and we need to
show the limit of $J_1$. We have  $\lambda f(u_1/\lambda, 0) = \lambda L(1,0) |u_1/\lambda| (1 + \delta(u_1/\lambda))
= L(1,0) |u_1| (1 + \delta(u_1/\lambda))$ where $\delta(u)$ is a bounded function tending to 0  as $u \to  0$. Therefore,
$J_1 =  J'_1 + J''_1 $, where
$$
J'_1 := \frac{2L(1,0)}{\log \lambda} {\rm Re} \int_{0}^{\lambda \pi}
(1- \e^{\i x_1 u})(1- \e^{-\i y_1 u}) u^{-1}  \d u \ \to \  2L(1,0) \begin{cases}2, &x_1 = y_1, \\
1, &x_1 \ne y_1
\end{cases}
$$
as  $\lim_{\lambda \to \infty} \int_1^\lambda  \e^{\i x u}  u^{-1} \d u$ exists for any $x \ne 0$,
and $|J''_1| \le C(\log \lambda)^{-1} `\int_0^{\lambda \pi} (u \wedge 1)^2 u^{-1} \delta (u/\lambda) \d u  \to  0 $ follows as in
\eqref{R2prime} by splitting the last integral over sets $u/\lambda < \epsilon $ and $\epsilon \le u/\lambda \le \pi $.
Whence, $R_{1}(\mbx, \mby) \to \kappa_+^2 \E  B_{0,1/2}(\mbx) B_{0,1/2}(\mby) $,  where
\begin{equation}\label{kappa5}
\kappa^2_+ = 4 L(1,0) \kappa^2 = 4 L(1,0) \int_{\R} |(1-\e^{\i u})/u|^2 \d u.
\end{equation}
To show that $\tilde R_1(\mbx, \mbx)$ is negligible,  use the decomposition of $\tilde f_1(\mbu)$ following \eqref{Jsub1} so that
$|\tilde R_1(\mbx, \mbx)| \le C\sum_{k=1}^3 \tilde R_{1,k}(\mbx, \mbx) $ where $\tilde R_{1,k}(\mbx, \mbx) = J_{\lambda,k}/\log \lambda $  and
$J_{\lambda,k} $ are as in the aforementioned proof. Then
$\tilde R_{1,1}(\mbx, \mbx) = o(1)$ as in \eqref{Jsub11} while 
$\tilde R_{1,2}(\mbx, \mbx) \le C \int_{\R} (1 + u_2^2)^{-1} \delta_\lambda (u_2) \d u_2 $, c.f. \eqref{Jsub12}, with
$$
\delta_\lambda (u_2) := \frac{1}{\log \lambda}
\int_{|u_1| \le \lambda \pi} (1 + |u_1|)^{-1} \big| L\big(\frac{u_1}{\lambda}, \frac{u_2}{\lambda^\gamma}\big) - L\big(\frac{u_1}{\lambda}, 0\big)\big|
\d u_1
$$
bounded and tending to 0 as $\lambda \to \infty  $ for any fixed $u_2 \ne 0$ by continuity of $\tilde L$  due to $\gamma \upsilon_2 > 1 $.
Finally, $\tilde R_{1,3}(\mbx, \mbx) \le  C(\log \lambda)^{-1} \int_{\lambda \Pi} (1  + |u_1|)^{-1} \d u_1 \int_{\R}   \delta(u_1/\lambda, u_2/\lambda^\gamma) (1 + u_2^2)^{-1} \d  u_2 +  o(1)  =  o(1) $  follows by splitting the last integral  over $|u_1| < \epsilon \lambda $ and
$ |u_1| > \epsilon \lambda $,  using $\int_{\epsilon \lambda }^{\pi \lambda} u_1^{-1} \d u_1 = \log (1/\epsilon) < \infty $ for any  small
$\epsilon >0$.

Next, let   $\gamma < 1/\upsilon_2, d^2_{\lambda, \gamma} = \lambda^{2H(\gamma)},
2H(\gamma) =  1 + \gamma(1- \upsilon_2)   $.  Note
$\lambda^{\gamma \upsilon_2} f(u_1/\lambda, u_2/\lambda^\gamma) \to f_0(0,u_2) =  L(0,1) |u_2|^{\upsilon_2}$. Then
$R_{\lambda, \gamma}(\mbx, \mby) \to \kappa^2_- \E B_{1/2, H^-_2}(\mbx) B_{1/2,H^-_2}(\mby)$ follows as in the case
$\upsilon_1 \wedge \upsilon_2 < 1, \gamma < \gamma_0 $,
with $\kappa^2_- = L(0,1) \int_{\R^2} |u_2|^{\upsilon_2} \prod_{j=1}^2 |(1-\e^{{\i} u_j})/u_j|^2  \d \mbu, $ 
c.f. \eqref{kappa3}.

 Let   $\gamma =  1/\upsilon_2,$ then $ d^2_{\lambda,\gamma} = \lambda^{2 H(\gamma)}
\log \lambda, \, 2H(\gamma) = \frac{1}{\upsilon_2}  $. We have
\begin{eqnarray}\label{Jsub10}
R_{\lambda, \gamma}(\mbx, \mby)
&\sim&\frac{1}{\log \lambda} \int_{\lambda \Pi \times \lambda^{1/\upsilon_2} \Pi} \prod_{i=1}^2 \frac{(1- \e^{\i x_i u_i})(1- \e^{-\i y_i u_i})}{|u_i|^2}
\lambda f(u_1/\lambda, u_2/\lambda^{1/\upsilon_2}) \d \mbu
\end{eqnarray}
where $\lambda f(u_1/\lambda, u_2/\lambda^{1/\upsilon_2}) = f_0(\mbu) (1 + \delta(u_1/\lambda, u_2/\lambda^{1/\upsilon_2})) $
  $= \sum_{k=1}^3 f_k(\mbu), $
$f_1(\mbu) := L(\mbu) |u_1|,  f_2(\mbu) := L(\mbu) |u_2|^{\upsilon_2},  $  $ f_3(\mbu) :=
 f_0(\mbu)\delta(u_1/\lambda, u_2/\lambda^{1/\upsilon_2}), $  and $\lim_{\mbu \to \0} \delta(\mbu) = 0$.  Accordingly,
$R_{\lambda, \gamma}(\mbx, \mby) \sim \sum_{k=1}^3  R_k(\mbx, \mby)$; we will show that
$R_1(\mbx, \mby)$ is the main term and $R_k(\mbx, \mby) \to 0, k=2,3$. Indeed,
\begin{eqnarray*}
R_1(\mbx, \mby)
&=&\frac{1}{\log \lambda} \int_{-\lambda \pi}^{\lambda \pi}
\frac{(1- \e^{\i x_1 u})(1- \e^{-\i y_1 u})}{|u|} h_\lambda(u) \d u, 
\end{eqnarray*}
where  $h_\lambda (u) := \int_{|w| \le \lambda^{1/\upsilon_2} \pi}
(1- \e^{\i x_2 w})(1- \e^{-\i y_2 w})|w|^{-2} L(u,w) \d w  \to h(u)\, (\lambda \to \infty) $ and where
\begin{eqnarray}\label{hu}
h(u)&:=& \int_{\R}
\frac{(1- \e^{\i x_2 w})(1- \e^{-\i y_2 w})}{|w|^2} L(u,w) \d w \\
&\to&L(1,0) \int_{\R} \frac{(1- \e^{\i x_2 w})(1- \e^{-\i y_2 w})}{|w|^2} \d w  = L(1,0) \kappa^2 (x_2 \wedge y_2) \nn
\end{eqnarray}
as $|u| \to  \infty $,  see \eqref{BHH}, \eqref{kappa} for the last equality. Whence, the convergence
$$
R'_1 (\mbx, \mby)\ :=  \ \frac{1}{\log \lambda} \int_{-\lambda \pi}^{\lambda \pi}
\frac{(1- \e^{\i x_1 u})(1- \e^{-\i y_1 u})}{|u|} h(u) \d u \  \to  \ \kappa_+^2  \E B_{0,1/2}(\mbx)  B_{0,1/2}(\mby)
$$
with $\kappa^2_+$ in \eqref{kappa5}
follows as in \eqref{Jsub6} and the same limit for $R_1 (\mbx, \mby)$ requires few changes.
Next, $|R_2(\mbx, \mbx)| \le  C (\log \lambda)^{-1}
\int_{\R^2} (1 + u_1^2)^{-1} (1 +  u_2^2)^{-1} |u_2|^{\upsilon_2} \d \mbu = O(1/\log \lambda)$. Finally,
$|R_3(\mbx, \mbx)| \le  C (\log \lambda)^{-1} \big( \int_0^{\lambda \pi} (1 + u)^{-1} \int_{\R} (1+ w^2)^{-1}
\delta (u/\lambda, w/\lambda^{1/\upsilon_2}) \d w  + O(1) \big) \to 0 $ follows as 
in the case  $\gamma > 1/\upsilon_2$.

\smallskip

\noi \underline {\it Case $\upsilon_1 =1, \upsilon_2 > 1, \gamma_0 = 1  $. }  Let $\gamma \ge  1, d^2_{\lambda, \gamma} =  \lambda^\gamma \log_+ \lambda$.
Then
\begin{eqnarray*}
R_{\lambda, \gamma}(\mbx, \mby)
&\sim&\frac{1}{\log \lambda}
\int_{\lambda \Pi \times \lambda^\gamma \Pi} \prod_{i=1}^2 \frac{(1- \e^{\i x_i u_i})(1- \e^{-\i y_i u_i})}{|u_i|^2}
\lambda f(u_1/\lambda, u_2/\lambda^\gamma) \d \mbu
\end{eqnarray*}
where $\lambda f(u_1/\lambda, u_2/\lambda^\gamma) \to  f_0(u_1,0) = L(1,0) |u|$ due to $\gamma \upsilon_2 > 1 $.
Then $R_{\lambda, \gamma}(\mbx, \mby) \to \kappa^2_+ \E B_{0,1/2}(\mbx) B_{0,1/2}(\mby)  $ 
as in \eqref{Jsub6} with $\kappa^2_+$ in \eqref{kappa5}.
Next, let $\gamma < 1$. Then $d_{\lambda,\gamma} = \lambda^{1/2}$ and
\begin{eqnarray*}
R_{\lambda, \gamma}(\mbx, \mby) 
&\sim&
\int_{-\pi\lambda}^{\pi \lambda}   \frac{(1- \e^{\i x_1 u})(1- \e^{-\i y_1 u})}{|u|^2} \d u \times
\int_\Pi   \frac{ (1- \e^{\i \lambda^\gamma x_2 w}) (1- \e^{-\i \lambda^\gamma y_2 w}) }
{|1- \e^{\i v}|^2}  f(u/\lambda, w) \d w \\
&\sim&\kappa^2_- \, \E  B_{1/2,0}(\mbx) B_{1/2,0}(\mby)
\end{eqnarray*}
as in \eqref{Jsub2} with $\kappa^2_- $ in \eqref{kappa4}.

\smallskip

\noi \underline{\it Case $\upsilon_1 =\upsilon_2  = \gamma_0 = 1  $. } The convergence in \eqref{Rlim} for $\gamma \ne 1$ leading to
limits $V_+ =\kappa_+ B_{0,1/2}$ and
$V_- = \kappa_- B_{1/2, 0} $ by symmetry follow as in the case
$\upsilon_1 =1, \upsilon_2 > 1$, with
with $\kappa^2_+$ in \eqref{kappa5} and $\kappa^2_- = 4 L(0,1) \kappa^2 $.  Let us prove that for  $\gamma = 1$ \eqref{Rlim} tends
to the sum of the latter limits leading to
\begin{eqnarray}\label{Jsub9}
V_0&=&\kappa_+ B_{0,1/2} + \kappa_- B_{1/2, 0},   \quad \text{with} \quad   d^2_{\lambda,\gamma} = \lambda \log_+ \lambda, \nn
\end{eqnarray}
where $B_{0,1/2} $ and $B_{1/2, 0}$ are mutually independent. Proceeding  as in \eqref{Jsub10} we see that
$f_\lambda (\mbu) := \lambda f(u_1/\lambda, u_2/\lambda)$  $ \to f_0 (\mbu) = L(\mbu) (|u_1| + |u_2|) $ and $R_{\lambda, \gamma}(\mbx, \mby)$ behaves asymptotically
as the sum of two terms
\begin{eqnarray*}
R_k(\mbx, \mby)
&:=&\frac{1}{\log \lambda}
\int_{(\lambda \Pi)^2} L(u_1, u_2) |u_k|   \prod_{j=1}^2 \frac{(1- \e^{\i x_j u_j})(1- \e^{-\i y_j u_j})}{|u_j|^2}
 \d \mbu, \quad k=1,2
\end{eqnarray*}
which tend to $\kappa^2_+ \E B_{0,1/2}(\mbx) B_{0,1/2}(\mby)  $ and $\kappa^2_- \E B_{1/2,0}(\mbx) B_{1/2,0}(\mby)$, respectively.
We also need to check that the term $R_3(\mbx, \mby)$ corresponding to $  f_\lambda (\mbu) - f_0 (\mbu) $ is negligible, viz.,
$ |R_3(\mbx, \mby)| \le C  (\log \lambda)^{-1}   $  $ \int_{(\lambda \Pi)^2} (|u_1| +  |u_2|) $  $ \delta(u_1/\lambda, u_2/\lambda) $  $
\prod_{j=1}^2 (1 + |u_j|^2)^{-1} \d \mbu \to 0$. We omit these details being similar to  \eqref{Jsub10}.
This ends the proof of Part II and   Theorem \ref{thmGaussND}, too. \hfill  $\Box$

\begin{remark} \label{remND} {\rm Following the terminology in time series \cite{gir2012},
the asymptotic constants $\kappa^2_\pm $ may be dubbed `long-range variances'.
It is notable that the only case when $\kappa^2_\pm $
depend on the spectral density outside of the origin are \eqref{kappa4}, \eqref{kappa44} corresponding to $H^\pm_i = 0$ or
spectrum ND under `edge effects',
see Remark \ref{remedge}. Obviously, these expressions for $\kappa^2_\pm $ make sense for continuous $f$ (the continuity can be relaxed)
but not for an arbitrary integrable or bounded $f$.  On the other hand, continuity of $f$ is a consequence of summability of covariance
function hence occurs under covariance SRD and ND.

}
\end{remark}

\section{Long-range negative and hyperbolic dependence}

Recall the asymptotic form  of the spectral density under Assumption  (F)$_{{\rm LRND},2}$:
\begin{equation}\label{FLRND}
f(\mbu)\  \sim  \  f_0(\mbu) =  \frac{L(\mbu)}{|u_1|^{\upsilon_1} + |u_2|^{\upsilon_2}}
\qquad |\mbu| \to 0,
\end{equation}
where  $L(\mbu) = L(-\mbu), \mbu \in  \R^2_+$ is a continuous generalized invariant function  such that 
\begin{equation}\label{LLRND}
L(\mbu) \sim \ell \Big(\frac{|u_2|}{\rho (\mbu)^{1/\upsilon_2}}\Big)^{\mu} \qquad (u_2 \to 0)
\end{equation}
for some $\mu \in (0,1), \ell >0 $. Define
\begin{eqnarray}
H^+_1&:=&\frac{1}{2}\Big(1 + \big((\upsilon_1 + \frac{\mu \upsilon_1}{\upsilon_2})  \wedge 1\big)\Big),
\quad H^+_2 :=   \frac{1}{2}\Big(1 -  \big(\mu  \wedge (\frac{\upsilon_2}{\upsilon_1}- \upsilon_2)\big)\Big),  \label{HLRND} \\
H^-_1 &:=&\frac{1}{2}
\big(1 + \upsilon_1 - \frac{\upsilon_1}{\upsilon_2 \vee 1}\big), \hskip2cm  H^-_2 := \frac{1}{2}\big(1 + (\upsilon_2 \wedge 1)\big)\nn
\end{eqnarray}
Note $H^-_i, i=1,2 $ in \eqref{HLRND} are the same as in Theorem  \ref{thmGaussLRD}, \eqref{H+},
whereas  the expressions for  $H^+_i, i=1,2 $  in \eqref{HLRND} and \eqref{H+} agree
if and only if $\mu = 0$. Also note that
$H^\pm_1,  H^-_2 \in [1/2,1] $ while   $H^+_2 = \frac{1}{2} (1+ \upsilon_2 - \frac{\upsilon_2}{\upsilon_1})
\in [1/2,1] $ for $\upsilon_1 \ge 1 $; for  $\upsilon_1 \le 1 $ we have
$H^+_2 \in (0, 1/2] $. Finally,
$H(\gamma) $ in \eqref{HgammaLRDN} is a continuous function of $\gamma, \upsilon_i, i=1,2, \mu, $ its value
$H(\gamma_0) = \frac{1}{2}(1 + \upsilon_1 + \frac{\upsilon_1}{\upsilon_2}) $ at $\gamma = \gamma_0 := \frac{\upsilon_1}{\upsilon_2} $
being the same as in Theorem  \ref{thmGaussLRD} independently of  $\mu$.

The following theorem excludes
some particular cases of parameters $\mu, \upsilon_i, i=1,2 $ which may require extra logarithmic
normalizing factor. It also leaves open the question about the scaling limits when both  Assumptions (F)$_{\rm LRN\!D, 1}$ and
(F)$_{\rm LRN\!D, 2}$ are satisfied.

\begin{theorem} \label{thmGaussLRND}
Let $X$ in \eqref{Sgamma} be a stationary linear RF on $\Z^2$ in \eqref{Xlin} with spectral density $f$ 
satisfying Assumption (F)$_{LRN\!D, 2}$.
In addition, let $\upsilon_1 \ne 1, \upsilon_1 + \frac{\upsilon_2}{\upsilon_1} \ne 1$.
Then:

\begin{itemize}

\item The scaling limits in \eqref{VgammaD} exist for any $\gamma >0$ and satisfy  \eqref{Vtrans}  with
$\gamma_0 = \frac{\upsilon_1}{\upsilon_2} $  and the unbalanced limits are given by
\begin{eqnarray}\label{VLRND}
V_+
&:=&\kappa_+ B_{H^+_1, H^+_2}, 
\qquad
V_-
\ :=\  \kappa_- B_{H^-_1, H^-_2},  
\end{eqnarray}
where $H^\pm_i, i=1,2$ as in \eqref{HLRND}.
The asymptotic constant $\kappa_+ $  is defined in \eqref{kappaLRND}, \eqref{kappaLRND1}, whereas $\kappa_- $ is the same as in
Theorem \ref{thmGaussLRD}.

\item  The well-balanced limit $V_0 := V_{0,1/\rho}  $ is given in \eqref{V0G}.

\item  The normalization in \eqref{VgammaD} is given by
$d_{\lambda, \gamma} :=    \lambda^{H(\gamma)}  $
with
\begin{eqnarray} \label{HgammaLRDN}
H(\gamma)&:=&
\begin{cases}  H^+_1 + \gamma H^+_2, &\gamma \ge \gamma_0,  \\
H^-_1 + \gamma H^-_2, &\gamma \le \gamma_0.
\end{cases}
\end{eqnarray}

\item The RF $X$ exhibits  scaling transition at $\gamma_0 = \frac{\upsilon_1}{\upsilon_2}$.

\end{itemize}
\end{theorem}

\medskip

\noi {\bf Proof.} Since \eqref{HgammaLRDN} satisfy $H(\gamma ) > 1/2 \,  (\forall \, \gamma >0)$,
the Lindeberg criterion in \eqref{LindX} is satisfied as in Theorem \ref{thmGaussLRD}.
For $\gamma \le \gamma_0$ the results of Theorem \ref{thmGaussLRND} and their proof completely agree with those
of Theorem \ref{thmGaussLRD}.  The subsequent proof on \eqref{Rlim}
is limited to  $\gamma > \gamma_0$ and split into three cases as follows.

\smallskip

\noi \underline{Case
$\upsilon_1 + \frac{\upsilon_1}{\upsilon_2} < 1 $.} 
According to \eqref{FLRND}-\eqref{LLRND} for
 $\gamma > \upsilon_1/\upsilon_2$ we have that
\begin{eqnarray*}
f(u_1/\lambda, u_2/\lambda^{\gamma}) 
&\sim&\frac{\lambda^{\upsilon_1}}{|u_1|^{\upsilon_1}} \times \ell
\Big|\frac{\frac{u_2}{\lambda^\gamma}}{|\frac{u_1}{\lambda}|^{\upsilon_1/\upsilon_2} } \Big|^{\mu} 
\ = \ \ell \lambda^{2H(\gamma) - 1 - \gamma}  \prod_{j=1}^2 |u_j|^{1-2H^+_j}
\end{eqnarray*}
point-wise in $\mbu = (u_1,u_2) \in \R^2_0  $. Then, the limit in \eqref{Rlim} with  $V_\gamma = V_+
= \kappa_+ B_{H^+_1, H^+_2} $ follows similarly as in Theorem \ref{thmGaussLRD} with 
\begin{eqnarray}
\kappa^2_+
&=&\ell  \prod_{j=1}^2  \int_{\R} |(1- \e^{\i u})/u|^2  |u|^{1 - 2 H^+_{j}}   \d u. 
\label{kappaLRND}
\end{eqnarray}

\medskip

\noi \underline{Case
$\upsilon_1 < 1  <\upsilon_1 + \frac{\upsilon_1}{\upsilon_2}$.} Then $(2 H^+_1, 2H^+_2)
= (1 +  \upsilon_1 + \frac{\mu \upsilon_1}{\upsilon_2}, 1- \mu) $ if $\mu < \frac{\upsilon_2}{\upsilon_1} - \upsilon_2$,
$= (2, 1 + \upsilon_2 - \frac{\upsilon_2}{\upsilon_1}) $ if $\mu >  \frac{\upsilon_2}{\upsilon_1} - \upsilon_2 $.
The proof of the limit $V_+
= \kappa_+ B_{H^+_1, H^+_2} $ for $\mu < \frac{\upsilon_2}{\upsilon_1} - \upsilon_2$ is  the same as
in the case $\upsilon_1 + \frac{\upsilon_1}{\upsilon_2} < 1 $ above, and is omitted. Next,
let $\mu >  \frac{\upsilon_2}{\upsilon_1} - \upsilon_2. $ Due to $\gamma > \upsilon_1/\upsilon_2$  we see that
$$
\lambda^{-\gamma \upsilon_2} f\big(\frac{u_1}{\lambda^{\gamma \upsilon_2/\upsilon_1}}, \frac{u_2}{\lambda^\gamma}\big) \to f_0(\mbu), \quad
\lambda^{-1} D_{[\lambda x]}\big(\frac{u_1}{\lambda^{\gamma \upsilon_2/\upsilon_1}}\big) \to x, \quad
\lambda^{-\gamma} D_{[\lambda^\gamma x]}\big(\frac{u_2}{\lambda^{\gamma}}\big) \to   \frac{1- \e^{\i x u_2}}{-\i u_2}
$$
point-wise for any $u_i \ne 0, i=1,2, x >0 $.  We also have $\int_{\R} g_+(u_1,u_2) \d u_1 =
|u_2|^{1 - 2 H^+_2} \int_{\R} f_0(u,1) \d u $ where the last integral converges due to \eqref{LLRND} and $\mu >  \frac{\upsilon_2}{\upsilon_1} - \upsilon_2$. Then,  $V_+
= \kappa_+ B_{1, H^+_2} $ follows as in \cite[Thm.3.1]{ps2015}. 
The asymptotic constant $\kappa_+ $ in both cases of $\mu$ is given by
\begin{eqnarray} \label{kappaLRND1}
\kappa^2_+
&=&
\begin{cases}
\ell  \prod_{j=1}^2  \int_{\R} |(1- \e^{\i u})/u|^2  |u|^{1 - 2 H^+_{j}}   \d u, &\mu  < \frac{\upsilon_2}{\upsilon_1} - \upsilon_2 , \\
\int_{\R} \big|(1- \e^{\i v})/v\big|^2 |v|^{1-2 H^+_2} \d v
\times \int_{\R} g_+ (u, 1) \d u, &\mu > \frac{\upsilon_2}{\upsilon_1} - \upsilon_2. 
\end{cases}
\end{eqnarray}

\vskip.1cm

\noi \underline{Case
$\upsilon_1  > 1 $.} In this case, the results do not depend on $\mu$ and completely agree with those in Theorem \ref{thmGaussLRD},
including the proof. This ends the proof of Theorem \ref{thmGaussLRND}.  \hfill $\Box$

\medskip

Next, we formalize the meaning of hyperbolic dependence mentioned in the Introduction.
Let 
\begin{equation} \label{fhyp}
f(\mbu) =  f_{0,{\rm  hyp}}(\mbu)(1 + o(1)), \quad |\mbu| \to 0, \quad \text{where} \quad
 f_{0,{\rm  hyp}}(\mbu)  :=  L(\mbu)\prod_{i=1}^2 |u_i|^{-\upsilon_i}, \qquad \mbu \in \Pi^2,
\end{equation}
$|\upsilon_i| <1, i=1,2 $ and
\begin{eqnarray}\label{Lhyp}
L(\mbu):= \tilde L \Big( \frac{u_1}{(|u_1|^{|\upsilon_1|} +  |u_2|^{|\upsilon_2|})^{1/|\upsilon_1|}},
\frac{u_2}{(|u_1|^{|\upsilon_1|} +  |u_2|^{|\upsilon_2|})^{1/|\upsilon_2|}} \Big)
\end{eqnarray}
is a generalized invariant  function corresponding to generalized homogeneous function
$\rho (\mbu) =  |u_1|^{|\upsilon_1|} +  |u_2|^{|\upsilon_2|}, \mbu = (u_1,u_2) \in \R^2_0 $.  Let
\begin{equation} \label{Hhyp}
H_i := \frac{1}{2}  (1 + \upsilon_i), \qquad i=1,2.
\end{equation}
The class in \eqref{fhyp} includes (separately) fractionally integrated spectral densities $\prod_{i=1}^2 |1- \e^{\i u_i}|^{-\upsilon_i} $
which play an important role in spatial statistics \cite{bois2005, guo2009, leo2013}.
If $L(\mbu) $ is separated from 0 and $\infty $, the spectral density in \eqref{fhyp} satisfies (spectrum) LRD, ND or LRND properties 
depending on the sign of $\upsilon_i, i=1,2 $ except that it explodes/vanishes on the coordinate axis $u_i =  0$ when
$\upsilon_i \ne 0$ as well and represents a different class from those discussed in Theorems \ref{thmGaussLRD}-\ref{thmGaussLRND}
Introduce a Gaussian RF
\begin{eqnarray} \label{Vhyp}
V_{\rm 0, hyp}(\mbx)&:=&\int_{\R^2} \prod_{j=1}^2 \frac{1-\e^{{\i} u_j x_j}}{\i u_j}\, \sqrt{f_{\rm 0, hyp}(\mbu)} Z(\d \mbu),
\quad \mbx \in \R^2_+,
\end{eqnarray}
where $Z(\d \mbu)$ is the same white noise as in \eqref{V0G}. In the case when $L(\mbu) = \ell >0$ is a constant
function, the RF  $V_{\rm 0, hyp}$ is a multiple of FBS: $V_{\rm 0, hyp} = \kappa \ell B_{H_1,H_2} $ with $H_i, i=1,2 $ given
in \eqref{Hhyp}  and $\kappa >0 $ as in \eqref{kappa}.

\begin{theorem}\label{farimaLM}
be a stationary linear RF on $\Z^2$ in \eqref{Xlin} with spectral density $f$  in  \eqref{fhyp}-\eqref{Lhyp},
where $\upsilon_i \in (-1,1), i=1,2$ and  $\tilde L(\mbx)$ is a strictly positive continuous function.
Then:
\begin{itemize}

\item The scaling limits in \eqref{VgammaD} exist for any $\gamma >0$ and satisfy 
\begin{equation}\label{Vgammahyp}
V_\gamma = \begin{cases}
\kappa_\gamma B_{H_1,H_2},  &\gamma \neq \frac{|\upsilon_1|}{|\upsilon_2|}, \\
V_{\rm 0, hyp}, &\gamma = \frac{|\upsilon_1|}{|\upsilon_2|},
\end{cases}
\end{equation}
where $H_i, i=1,2$ as in \eqref{Hhyp} and
\begin{equation}\label{kappahyp}
\kappa^2_\gamma := \prod_{j=1}^2 \int_{\R} |(1- \e^{\i u})/u|^2 |u|^{1-2H_j} \d u \times
\begin{cases}
L(1,0), &\gamma  > \frac{|\upsilon_1|}{|\upsilon_2|}, \\
L(0,1), &\gamma  < \frac{|\upsilon_1|}{|\upsilon_2|}.
\end{cases}
\end{equation}

\item  The normalization in \eqref{VgammaD} is given by
$d_{\lambda, \gamma} :=    \lambda^{H(\gamma)}, H(\gamma)  := H_1 + \gamma H_2.  $

\item The RF $X$  does not exhibit  scaling transition.

\end{itemize}

\end{theorem}

We omit the proof of Theorem \ref{farimaLM} since it resembles \cite{ps2015} and
the  previous proofs. The gaussianity criterion \eqref{LindX} holds for ${\rm sign}(\upsilon_1) = {\rm sign}(\upsilon_2) $ as
in Theorems \ref{thmGaussLRD} and \ref{thmGaussND}; for $\upsilon_1 \upsilon_2 < 0$ from \eqref{fhyp} we get
$|\widehat a(\mbu)| \le C \prod_{i=1}^2 |u_i|^{-\upsilon_i/2} $ and   \eqref{LindX} follows similarly.
The absence of scaling transition at $\gamma = \frac{|\upsilon_1|}{|\upsilon_2|}$ is clear from
\eqref{Vgammahyp}-\eqref{kappahyp} and the fact that $L(1,0) >0, L(0,1) >0$ according to  the positivity assumption on $L$ 
(the last conclusion may fail if $L(1,0)$ and/or $L(0,1)$ vanish).

\smallskip

We end the paper with two remarks on possible extensions of this work.

\begin{remark} \label{remedge}  {\rm  {\it More general scaling schemes and `edge effects'.}
In a more abstract setting, a RF is often indexed by test functions, through which scaling operations are defined; see
\cite{dob1979, dob1980, gelf1964}. Accordingly, one can study anisotropic scaling limits of integrals 
\begin{equation}\label{Sphi}
S_{\lambda,\gamma}(\phi)  :=  \int_{\R^2} \phi (\lambda^{-\Gamma} \mbt)
X(\lceil \mbt \rceil ) \d \mbt
\end{equation}
 involving $X $ in \eqref{Xlin} extended to $\R^2$; $ \lceil \mbt \rceil := (\lceil t_1 \rceil,
 \lceil t_2 \rceil), \mbt = (t_1,t_2) \in \R^2 $,
 and a re-scaled function $\phi: \R^2 \to  \R$, $\lambda^{-\Gamma} \mbt :=  (\lambda^{-1} t_1, \lambda^{-\gamma}  t_2), \Gamma := {\rm diag}(1, \gamma)$, for a given $\gamma >0$ and each $\phi $ from a class $\Phi$ of (test) functions. The sum  in \eqref{Sgamma} corresponds to \eqref{Sphi}
with indicator function $\phi (\mbt) = \I(\mbt \in ]\boldsymbol{0}, \mbx])$. For \eqref{Sphi},
the scaling  limit $d^{-1}_{\lambda, \gamma} (S_{\lambda,\gamma}(\phi) - \E S_{\lambda,\gamma}(\phi)) \limd V_\gamma(\phi) $
is a RF indexed by  $\phi \in \Phi$.  Isotropic $(\gamma = 1)$ scaling limits in spirit of \eqref{Sphi} for different classes
of RF on $\R^d $ were studied in \cite{dob1980, bier2010, kaj2007} and other papers. Extending
Theorems \ref{thmGaussLRD}-\ref{farimaLM} to scaling limits of integral functionals in \eqref{Sphi} for suitable class $\Phi$ of test functions
which include indicator functions $\phi(\mbt) = \I(\mbt \in A)$  of general bounded sets $A \subset \R^2 $
is an interesting problem.
Of particular interest is the case of ND RF $X$,  where `edge effects' can be expected following \cite{lah2016, sur2020},
leading to scaling limits $V_\gamma $  `living'
on the boundary of $A$ or the discontinuity set of $\phi$.

}
\end{remark}

\begin{remark} \label{remincong} {\rm {\it Incongruous scaling and dependence axis.} For linear LRD RF  $X$ in \eqref{Xlin} \cite{pils2020}
define the {\it dependence axis of $X$} as the direction on the plane 
along which $a(\mbt)$ decay at the smallest rate. 
In spectral  terms, we may define the dependence axis as a direction along which the spectral density $f(\mbx)$
grows at the fastest  rate when $|\mbx| \to 0  $.
Under Assumption (F)$_{\rm  LRD}$ with $\upsilon_1 \ne \upsilon_2 $  such dependence axis 
coincides with one of the coordinate axes.  
The scaling in \eqref{VgammaD} involves
rectangles with sides parallel to the
coordinate  axis.
Using the terminology  in  \cite{pils2020} we may say
that the scaling in \eqref{VgammaD} and in
Theorem  \ref{thmGaussLRD} is   {\it congruous with the dependence axis of $X$}. 
\cite{pils2020} showed that {\it incongruous scaling} of linear LRD RF  in \eqref{Xlin}  
may dramatically change the scaling limits in  \eqref{VgammaD}  and the scaling transition point $\gamma_0$. 
We expect that the spectrum approach in our paper may lead to 
a comprehensive treatment of incongruous scaling limits 
under various dependence assumptions.  

}
\end{remark}

\section*{Acknowledgements}

The author thanks Vytaut\.e Pilipauskait\.e for useful comments and Remigijus Lapinskas
for help with Figure 1 graphs.

\end{document}